\documentclass[english,12pt]{article}
\usepackage[T1]{fontenc}
\usepackage[utf8]{inputenc}

\usepackage{fullpage}
\usepackage{amsmath,amsthm}

\newtheorem{thm}{Theorem}
\numberwithin{thm}{section}
\newtheorem{lem}{Lemma}
\numberwithin{lem}{section}

\theoremstyle{definition}
\newtheorem{defi}{Definition}
\numberwithin{defi}{section}

\numberwithin{exper}{section}
\theoremstyle{exper}

\theoremstyle{definition}
\newtheorem{assumption}{Assumption}
\numberwithin{assumption}{section}

\theoremstyle{definition}

\numberwithin{exam}{section}

\theoremstyle{definition}
\newtheorem{remark}{Remark}
\numberwithin{remark}{section}

\numberwithin{hequa}{section}

\usepackage{amssymb}
\usepackage{graphicx}
\usepackage{subfigure}
\usepackage{color}
\usepackage{epsfig}
\usepackage{epstopdf}
\usepackage{hyperref}
\usepackage{algorithm,algcompatible}
\usepackage{authblk}
\usepackage{comment}
\usepackage{footmisc} 
\usepackage{hyperref} 

\numberwithin{equation}{section}

\graphicspath{ {./figures/} }

\title{Order Reduction of Exponential Runge--Kutta Methods: Non-Commuting Operators}

\author{Trung-Hau Hoang\thanks{\href{mailto:trunghaugg@gmail.com}{trunghaugg@gmail.com}}}
\affil{Department of Mathematics, University of Innsbruck, Austria} 

\date{ }

\begin{document}

\maketitle

\section*{Abstract}

Nonlinear parabolic equations are central to numerous applications in science and engineering, posing significant challenges for analytical solutions and necessitating efficient numerical methods. Exponential integrators have recently gained attention for handling stiff differential equations. This paper explores exponential Runge--Kutta methods for solving such equations, focusing on the simplified form $u^{\prime}(t)+A u(t)=B u(t)$, where $A$ generates an analytic semigroup and $B$ is relatively bounded with respect to $A$. By treating $A$ exactly and $B$ explicitly, we derive error bounds for exponential Runge--Kutta methods up to third order. Our analysis shows that these methods maintain their order under mild regularity conditions on the initial data $u_0$, while also addressing the phenomenon of order reduction in higher-order methods. Through a careful convergence analysis and numerical investigations, this study provides a comprehensive understanding of the applicability and limitations of exponential Runge--Kutta methods in solving linear parabolic equations involving two unbounded and non-commuting operators.

\section{Introduction}

Numerous applications in science and engineering can be modeled using the following nonlinear parabolic equation for the unknown $u(t,x)$:
\begin{equation}\label{eq98} 
	\begin{aligned}
		& \partial_{t}u - \triangle u = f(  \nabla u, u), \qquad u(0) = u_0, 
	\end{aligned}	
\end{equation}
where \( f \) is a nonlinear function of \( \nabla u \) and \( u \). Due to the complexity of the problem, an analytical solution to \eqref{eq98} is usually not available. In this case, numerical methods play a crucial role in approximating the solution to \eqref{eq98}. Recently, exponential integrators (\cite{HO2010}) have become an attractive numerical method for efficiently integrating systems of differential equations that exhibit stiffness (see \cite{10.1007/s10543-013-0446-0, CROUSEILLES2020109688, DEKA2023101302}). These methods solve the linear part exactly and integrate the nonlinearity explicitly. 

When $f$ depends only on $u$, error analysis of exponential Runge--Kutta methods has been derived (\cite{HO2005}). In this paper, we focus on the study of exponential Runge--Kutta methods, which apply to a simplified form of \eqref{eq98} as follows
\begin{equation}\label{4121}
	u^{\prime}(t) + Au(t) = Bu(t), \quad u\left(0\right)=u_0.    
\end{equation}
Here, \( A \) is an operator generating an analytic semigroup, and \( B \) is an operator relatively bounded with respect to \( A \). For example, \( A \) could be a second-order strongly elliptic operator, and \( B \) could be a first-order differential operator. Studying the linear problem \eqref{4121} is a fundamental step towards understanding the behavior of the nonlinear problem \eqref{eq98}.

The exact solution of \eqref{4121} can be written as $u(t) = \mathrm{e}^{-t(A-B)}u_0$. However, we choose to treat the dominant part $A$ exactly, whereas treating the unbounded part $B$ explicitly.~For example, using the \textit{variation-of-constants} formula and integrating the nonlinearity at the only known value $u_0$, yields the exponential Euler method
\begin{equation*}
	u_{1}=\mathrm{e}^{-\tau A} u_0+\tau \varphi_1(-\tau A) Bu_0, \qquad \varphi_1(x) = \frac{e^x-1}{x}.	
\end{equation*}
There are two main reasons for this approach. First, we aim to investigate the order of exponential Runge--Kutta methods for the unbounded operator $B$. Second, computing  $u(t) = \mathrm{e}^{-t(A-B)}u_0$ can be expensive due to the complex structure of $A-B$ (see \cite{doi:10.1137/23M1562056}), while approximating $\mathrm{e}^{-tA}u_0$ is more straightforward. In this paper, we are interested in the exponential Runge--Kutta methods up to the third order (see \cite{HO2005, HO2010}). 

This paper is structured as follows. The convergence analysis of \eqref{4121} is carried out within the standard framework of analytic semigroups in a Banach space \( X \), which is recalled in Section \ref{Analyticalframework}. We briefly review explicit exponential Runge--Kutta methods applied to \eqref{4121} in Section \ref{EXPRK}. In Section \ref{erroranalysis} we analyse the convergence of explicit exponential Runge--Kutta methods applied to \eqref{4121}. The convergence results for first-order and second-order methods are given in Theorems \ref{theo2} and \ref{theo3} respectively, with the primary result detailed in Theorem \ref{theo4}. Section \ref{numericalchap2} presents numerical investigations demonstrating possible theoretical order reductions. Finally, conclusions are contained in Section \ref{concluchap2}.

 \section{Analytical framework}\label{Analyticalframework} 

In this section, we present an analytical framework that serves as the basis for the forthcoming convergence analysis. Let $X$ be a Banach space, and let $\mathcal{D}(A)$ denote the domain of $A$ in $X$. Our assumptions on the operators $A$ and $B$ are those of \cite{henry1981geometric,pazy1983semigroups}.

\begin{assumption}\label{ass1}
	Let $A: \mathcal{D}(A) \rightarrow X$ be sectorial; i.e., $A$ is a closed linear operator in $X$ and densely defined such that, for $\phi \in $  $(0, \pi / 2)$, $M \in [1, \infty)$ and $a\in \mathbb{R}$, the following conditions holds
		$$
		\left\|(\lambda I-A)^{-1}\right\| \leq \frac{{M}}{|\lambda-{a}|}  \quad \text { for all } \quad \lambda \in \mathrm{S}_{{a}, \phi} ,
		$$
		for the sector
		$$\qquad \mathrm{S}_{{a}, \phi}=\{\lambda \in \mathbb{C} \text{ }|\text{ }\phi \leq| \arg (\lambda-{a})|  \leq \pi,\text{ } \lambda \neq {a}\} $$
		that lies in the resolvent set of $A$.  
\end{assumption}

The operator $-A$ is an infinitesimal generator of an analytic semigroup $\left\{\mathrm{e}^{-t A}\right\}_{t \geq 0}$ under Assumption \ref{ass1}. For $\omega \in (-a, \infty)$, the fractional powers of $\widetilde{A} = A + \omega I$ are well-defined (see \cite{henry1981geometric}). For notational simplicity, we set $\omega = 0$, and hence $\widetilde{A} = A$. The following stability bounds are useful in the later error analysis (see \cite{HOCHBRUCK2005323}).

\begin{lem}\label{lem1}
	For fixed $\omega \in (-a,\infty)$ and together with  Assumption \ref{ass1}, the following bounds 
	\begin{equation}\label{parabolicsmoothing1}
			\left\|\mathrm{e}^{-t A}\right\|+\left\|t^\gamma A^\gamma \mathrm{e}^{-t A}\right\| \leq C, \quad \gamma \geq 0, 
	\end{equation}
	hold uniformly on $t \in [0,T]$.
\end{lem}

\begin{assumption}\label{ass2} 
Let \( 0 < \gamma \leq 1 \), and let \( B: \mathcal{D}(B) \rightarrow X \) be a closed linear operator satisfying
	\begin{equation}\label{eq99}
		\mathcal{D}(A^{\gamma}) \subset \mathcal{D}(B) \text{ and } \|B x\| \leq \varepsilon\|A^{\gamma} x\|+K(\varepsilon)\|x\| \text{ for all } x \in D(A^\gamma)
	\end{equation}
for sufficiently small $\varepsilon>0$. Here, \( K(\varepsilon) \) denotes a positive function. Under Assumption \ref{ass2}, there exists $C > 0$ such that (see 	\cite{pazy1983semigroups})   
\begin{equation}\label{eq14}
	\left\|  B A^{-\gamma} \right\| \leq C.
\end{equation}

In addition, we suppose the following 
\begin{equation}\label{eq15}
	\left\| A^{-\gamma} B  \right\| \leq C.
\end{equation}

\end{assumption}

\begin{defi}
	Let $M \in \mathbb{R}$ be an arbitrary number. We define $$M^- = M - \widehat{\zeta}  \text{ and } M^+ = M + \widetilde{\zeta} ,$$ where $\widehat{\zeta}, \widetilde{\zeta} > 0$ are any fixed small numbers. 
\end{defi}

\section{Exponential Runge--Kutta methods}\label{EXPRK}

We consider a class of explicit one-step methods known as the exponential Runge--Kutta methods (see \cite{HO2005}). These methods when applied to problem \eqref{4121} take the following form: 
\begin{subequations}\label{eq20}
	\begin{equation}\label{eq20a}
		U_{n i}  =\mathrm{e}^{-c_i \tau A}u_n +\tau \sum_{j=1}^{i-1} a_{i j}\left(-\tau A\right) B  U_{n j}, \quad 1 \leq i \leq s,
	\end{equation}
	\begin{equation}\label{eq20b}
		u_{n+1} =\mathrm{e}^{-\tau A }u_n +\tau\sum_{i=1}^s b_i\left(-\tau A\right) B U_{n i},
	\end{equation}
\end{subequations}
where $u_n\approx u(t_n)$ and $U_{n i} \approx u(t_n+c_i \tau)$. Here $t_{n+1} = t_{n} + \tau$, where $t_n = n\tau$, and $\tau > 0$ denotes a positive time step size. The coefficients $a_{i j}\left(-\tau A\right)$ and $b_i\left(-\tau A\right)$ are usually the linear combinations of the functions $\varphi_k(-\tau A)$. These functions take the form of
	$$
\varphi_k(z)=\int_0^1 \mathrm{e}^{(1-\theta) z} \frac{\theta^{k-1}}{(k-1) !} d \theta, \quad k \geq 1 ,
$$
and they fulfill the following recursion relation
$$
\qquad \varphi_0(z)=\mathrm{e}^z, \qquad \varphi_{k+1}(z)=\frac{\varphi_k(z)-\varphi_k(0)}{z}, \qquad k \geq 0 . 
$$
Using Lemma \ref{parabolicsmoothing1}, the functions $\varphi_k(-\tau A)$ are well-defined. The coefficients \( a_{ij}(-\tau A) \) and \( b_i(-\tau A) \) (denoted generically as \( \phi(-t A) \)) satisfy 
\begin{equation}\label{eq19}
	\|\phi(-t A)\| + \left\|t^\eta {A}^\eta \phi(-t A)\right\| \leq C, \quad 0 \leq \eta \leq 1.	
\end{equation}

The order conditions for the third order exponential Runge--Kutta methods have been derived (see \cite{HO2005}). These conditions are recalled in the Table \ref{Table1}.

\def\arraystretch{1.4}

\begin{table}[H]
	\begin{center}
		\begin{tabular}{|c|c|c|}
			\hline No. & Order &  Condition \\
			\hline 1 & 1 & $\sum_{i=1}^s b_i(Z)=\varphi_1(Z)$ \\
			\hline 2 & 2 & $\sum_{i=2}^s b_i(Z) c_i=\varphi_2(Z)$ \\
			3 & 2 & $\sum_{j=1}^{i-1} a_{i j}(Z)=c_i \varphi_1\left(c_i Z\right), \quad i=2, \ldots, s$ \\
			\hline 4 & 3 & $\sum_{i=2}^s b_i(Z) \frac{c_i^2}{2 !}=\varphi_3(Z)$ \\
			5 & 3 & $\sum_{i=2}^s b_i(Z) J \left( \sum_{k=2}^{i-1} a_{i k}(Z) c_k -c_i^2 \varphi_2\left(c_i Z\right) \right) =0$ \\
			\hline
		\end{tabular}	
		\caption{Stiff order conditions for explicit exponential Runge--Kutta methods up to order three. The variables $J$ and $Z$ denote arbitrary operators on $X$ (see \cite{HO2005}).  }
		\label{Table1}
	\end{center}
\end{table}
\section{Convergence results for exponential Runge--Kutta methods}\label{erroranalysis}

This section is dedicated to deriving error bound for the exponential Runge--Kutta methods \eqref{eq20} applied to \eqref{4121}. First, we derive the error recursion.

\subsection{Error recursion and defects representation}\label{errorrecursion}
The error recursion can be determined by inserting the exact solution into the numerical scheme to identify the defects in the stages, cf. \cite{HO2005, HOCHBRUCK2005323}. With the unbounded operator $B$, many complex new terms are formed, as we will see shortly.
 
  The exact solution of the initial value problem $\eqref{4121}$ can be represented using the variation-of-constants formula 
\begin{equation}\label{eq94} 
	u\left(t_n+\theta \tau \right)=\mathrm{e}^{-\theta \tau A} u\left(t_n\right)+\int_0^{\theta \tau} \mathrm{e}^{-(\theta \tau-\xi) A} Bu\left(t_n+\xi\right) \mathrm{d} \xi .
\end{equation}
By expanding $u$ in a Taylor series and considering the remainder part as an integral expression, we obtain
\begin{equation}\label{eq97}  
	u\left(t_n+\tau\right)=\sum_{j=1}^q \frac{\tau^{j-1}}{(j-1) !} u^{(j-1)}\left(t_n\right)+\int_0^\tau \frac{(\tau-\sigma)^{q-1}}{(q-1) !} u^{(q)}\left(t_n+\sigma\right) \mathrm{d} \sigma .
\end{equation} 
Substituting \eqref{eq97} into the right-hand-side of \eqref{eq94} gives
\begin{equation}\label{eq87}
	\begin{aligned}
		u\left(t_n+c_i \tau \right)= & \mathrm{e}^{-c_i \tau A} u\left(t_n\right)+\sum_{j=1}^{q_i}\left(c_i \tau \right)^j \varphi_j\left(-c_i \tau A\right) Bu^{(j-1)}\left(t_n\right) \\
		& +\int_0^{c_i \tau} \mathrm{e}^{-\left(c_i \tau-\xi\right) A} \int_0^\xi \frac{(\tau-\sigma)^{q_i-1}}{\left(q_i-1\right) !} Bu^{\left(q_i\right)}\left(t_n+\sigma\right) \mathrm{d} \sigma \mathrm{d} \xi .
	\end{aligned}
\end{equation}

Placing the exact solution in the numerical scheme \eqref{eq20} yields
\begin{subequations}\label{eq84}
	\begin{equation}\label{eq84a}
		u\left(t_n+c_i \tau \right)=\mathrm{e}^{-c_i \tau A} u\left(t_n\right)+\tau \sum_{j=1}^{i-1} a_{i j}(-\tau A) Bu\left(t_n+c_j \tau \right)+\Delta_{n i},
	\end{equation}
	\begin{equation}\label{eq84b}
		u\left(t_{n+1}\right)=\mathrm{e}^{-\tau A} u\left(t_n\right)+\tau \sum_{i=1}^s b_i(-\tau A) Bu\left(t_n+c_i \tau\right)+\delta_{n+1},
	\end{equation}
\end{subequations}
with defects $\Delta_{n i}$ and $\delta_{n+1}$. Inserting \eqref{eq97} in \eqref{eq84a}, we get
\begin{equation}\label{eq86}
	\begin{aligned}
		u\left(t_n+c_i \tau \right)= & \mathrm{e}^{-c_i \tau A} u\left(t_n\right)+ \tau \sum_{k=1}^{i-1} a_{i k}(-\tau A) \sum_{j=1}^{q_i} \frac{\left(c_k \tau \right)^{j-1}}{(j-1) !} Bu^{(j-1)}\left(t_n\right) \\
		& +\tau \sum_{k=1}^{i-1} a_{i k}(-h A) \int_0^{c_k \tau} \frac{\left(c_k \tau-\sigma\right)^{q_i-1}}{\left(q_i-1\right) !} Bu^{\left(q_i\right)}\left(t_n+\sigma\right) \mathrm{d} \sigma+\Delta_{n i} .
	\end{aligned}
\end{equation}

The explicit expression for the defects $\Delta_{n i}$ can be obtained by subtracting \eqref{eq87} from \eqref{eq86}. However, due to the presence of the unbounded operator $B$, many complex new terms arise, making them difficult to express. We will discuss each of these terms in the following subsections. The defects $\delta_{n+1}$ can be obtained in a similar way. 

Let $e_n=u_n-u\left(t_n\right)$ and $E_{n i}=U_{n i}-u\left(t_n+c_i \tau\right)$ denote the errors between the numerical solution and the exact solution. Subtracting \eqref{eq84} from the numerical method \eqref{eq20} leads to the following error recursion 
\begin{subequations}\label{eq85}
	\begin{equation}\label{eq85a}
		\quad E_{n i}=\mathrm{e}^{-c_i \tau A} e_n+\tau \sum_{j=1}^{i-1} a_{i j}(-\tau A)B E_{n j}-\Delta_{n i},
	\end{equation}
	\begin{equation}\label{eq85b}
		e_{n+1}=\mathrm{e}^{-\tau A} e_n+\tau \sum_{i=1}^s b_i(-\tau A)B E_{n i}-\delta_{n+1}.
	\end{equation}
\end{subequations}

 The error in the internal stages, \(E_{n i}\), is then entered into \eqref{eq85b} to find the global error expression.

 \subsection{Error bound for the exponential Euler method}\label{analysefirstorder} 
We begin with the numerical scheme of the simplest form of exponential Runge–Kutta methods: the exponential form of Euler's method (see \cite{HO2005}).~Applied to \eqref{4121}, this is as follows
\begin{equation}\label{Euler}
	u_{n+1}=\mathrm{e}^{-\tau A} u_n+\tau \varphi_1(-\tau A) Bu_n.
\end{equation}

For $s=1$, the following error recursion applies to the exponential Euler method

\begin{equation}\label{eq81}
	e_{n+1}=\mathrm{e}^{-\tau A} e_n+\tau \varphi_1(-\tau A) Be_n -\delta_{n+1}
\end{equation} 
where
 $$\delta_{n+1} = \int_0^{ \tau} \mathrm{e}^{-( \tau-\xi) A} \int_0^\xi Bu^{\prime}\left(t_n+\sigma\right) \mathrm{d} \sigma   \mathrm{d} \xi.  $$
  Solving error recursion \eqref{eq81}, we obtain
\begin{equation}\label{eq95}
	e_n= \tau\sum_{j=0}^{n-1} \mathrm{e}^{-(n-j-1) \tau A} \varphi_1(-\tau A) Be_j - \sum_{j=0}^{n-1} \mathrm{e}^{-j \tau A} \delta_{n-j} .		
\end{equation}
To bound $e_n$, we first establish the following lemma.
\begin{lem}\label{boundeduprime}
	Let Assumptions \ref{ass1} and \ref{ass2} be satisfied. Further, we assume that $u_0\in \mathcal{D}(A)$; then
	$$ \sup_{0 \leq t \leq T} \left\|  u^{\prime}(t) \right\|  \leq C .
	$$ 
\end{lem}

\begin{proof}
	The boundedness is a straightforward consequence of the estimate
	\begin{equation*}
		\left\| u^{\prime}(t) \right\|  \leq \left\| \mathrm{e}^{-t(A-B)}(A-B)u_0 \right\| \leq \left\| \mathrm{e}^{-t(A-B)} \right\| \left\| (A-B)u_0 \right\| \leq C,  
	\end{equation*}
$\text{for all } t \in [0,T]$. Henceforth, we note that \( A - B \) forms a new sectorial operator. The bounds for \(\mathrm{e}^{-t(A-B)}\) and \((A-B)u_0\) can be obtained from Assumptions \ref{ass1} and \ref{ass2}.
\end{proof}

The main idea in order to estimate all the terms is to distribute $A^{\gamma}$ and $A^{-\gamma}$ in an appropriate way. $A^{\gamma}$ usually appears next to the exponential term, and thanks to the commutativity of $\phi(-\tau A)$ with $A$, we can move $A^{-\gamma}$ next to the operator $B$. Using Assumption \ref{ass2}, the unbounded operator $B$ can be controlled. In addition, the exponential term together with $A$ can also be estimated using \eqref{parabolicsmoothing1} and \eqref{eq19}. For example, applying \eqref{parabolicsmoothing1} and \eqref{eq15}, the first part of \eqref{eq95} can be estimated by
\begin{equation*}
	\begin{aligned}
		&\left \| \tau\sum_{j=0}^{n-1} \mathrm{e}^{-(n-j-1) \tau A} A^{\gamma} \varphi_1(-\tau A) A^{-\gamma} Be_j  \right \| 
		\\
		& \leq  \tau\sum_{j=0}^{n-2} \left \| \mathrm{e}^{-(n-j-1) \tau A} A^{\gamma} \right \| \left \|  \varphi_1(-\tau A) \right \|  \left \| A^{-\gamma}  B \right \| \left \| e_j \right \| + \tau \left \| \varphi_1(-\tau A) A^{\gamma}   \right \| \left \|  A^{-\gamma} B      \right \| \left \|   e_{n-1}     \right \| \\
		& \leq C \tau  \sum_{j=0}^{n-2} t_{n-j-1}^{-\gamma}   \left \| e_j \right \| +  C \tau^{1 - \gamma} \left \| e_{n-1}     \right \| . 
	\end{aligned}
\end{equation*}
To estimate the second part of \eqref{eq95}, in addition to the techniques used to estimate the first part of \eqref{eq95}, Lemma \ref{boundeduprime} is also employed.~The term for $j=0$ can be bounded as
\begin{equation}\label{eq11} 
	\begin{aligned}
		 \left \| \delta_{n}  \right \|  & \leq  \int_0^{ \tau} \int_0^\xi \left \| \mathrm{e}^{-( \tau-\xi) A} A^{\gamma} \right \| \left \| A^{-\gamma} B \right \| \left \| u^{\prime} \left(t_{n-1}+\sigma\right)\right \|  \mathrm{d} \sigma   \mathrm{d} \xi   \\ 
		 & \leq C \sup_{0 \leq t \leq T} \left\|  u^{\prime}(t) \right\| \int_0^{ \tau} \int_0^\xi (\tau-\xi)^{-\gamma}    \mathrm{d} \sigma   \mathrm{d} \xi \leq C  \tau^{2 - \gamma}  
	\end{aligned}
\end{equation}
whereas the remaining sum is bounded by 
{\allowdisplaybreaks
	\begin{equation}\label{eq13}
	\begin{aligned}
		\left \| \sum_{j=1}^{n-1} \mathrm{e}^{-j \tau A} \delta_{n-j}  \right \|  
		& \leq \sum_{j=1}^{n-1} \int_0^{ \tau} \int_0^\xi \left \| \mathrm{e}^{-j\tau A} A^{\gamma} \right \| \left \| \mathrm{e}^{-( \tau-\xi) A} \right \| \left \| A^{-\gamma} B \right \| \left \| u^{\prime} \left(t_{n-j-1}+\sigma\right)\right \|  \mathrm{d} \sigma   \mathrm{d} \xi   \\   
		 & \leq C \sup_{0 \leq t \leq T} \left\|  u^{\prime}(t) \right\| \tau^2   \sum_{j=1}^{n-1} t_j^{-\gamma} 
		 \leq C  \tau    . 	
	\end{aligned}	
	\end{equation}
}At this point, the error bound for $e_n$ can be given.  However, the error bounds without assuming \( u_0 \in \mathcal{D}(A) \) can also be derived.~We establish the following lemma to support this.
\begin{lem}\label{lem7} 
	Let Assumptions \ref{ass1} and \ref{ass2} be satisfied, then
	$$  \left\|  u^{\prime}(t) \right\|  \leq Ct^{-1} .
	$$ 
\end{lem}
\begin{proof}
	The bound is obtained by the following estimation 
$$
		\left\| u^{\prime}(t) \right\|  \leq \left\| \mathrm{e}^{-t(A-B)}(A-B)u_0 \right\| \leq \left\| \mathrm{e}^{-t(A-B)}  (A-B)\right\| \left\| u_0 \right\| \leq Ct^{-1}.
$$ 
Here, to bound \(\mathrm{e}^{-t(A-B)} (A-B)\), we use a similar property as in \eqref{parabolicsmoothing1}.
\end{proof}
Using \eqref{parabolicsmoothing1}, \eqref{eq15}, and Lemma \ref{lem7}, the bound for the term with $j = n-1$ is given by 
\begin{equation*}
	\begin{aligned}
		\left \| \mathrm{e}^{-(n-1) \tau A} \delta_{1}  \right \| 
		& \leq \left \|  \mathrm{e}^{-(n-1) \tau A} A^{\gamma} \right \| \int_0^{ \tau} \int_0^\xi \left \| \mathrm{e}^{-( \tau-\xi) A}  \right \| \left \| A^{-\gamma} B \right \| \left \| u^{\prime} \left(\sigma\right)\right \|  \mathrm{d} \sigma   \mathrm{d} \xi  
		\\
		& \leq C t_{n-1}^{-\gamma} \int_0^{ \tau} \int_0^\xi \sigma^{-1}    \mathrm{d} \sigma   \mathrm{d} \xi .  
	\end{aligned}
\end{equation*}
Without assuming \( u_0 \in \mathcal{D}(A) \), we encounter a problem here because the integral $\int_0^\xi \sigma^{-1} \, \mathrm{d} \sigma$ is undefined.
Therefore, we have to derive a new formula for $\delta_{1}$.~Placing the exact solution in the numerical scheme \eqref{Euler} gives 
\begin{equation}\label{eq91}
	u(t_{1})=\mathrm{e}^{-\tau A} u_0+\tau \varphi_1(-\tau A) Bu_0+ \delta_{1}.
\end{equation}
By using the variation-of-constants formula, we also have 
\begin{equation}\label{eq92}
	u( t_1)=\mathrm{e}^{- \tau A} u_0+\int_0^{ \tau} \mathrm{e}^{-( \tau-\xi) A} Bu\left( \xi \right)  \mathrm{d} \xi.
\end{equation}
The new formulation for $\delta_{1}$ is obtained by subtracting \eqref{eq91} from \eqref{eq92}
\begin{equation}\label{eq93}
	\delta_{1} = \int_0^{ \tau} \mathrm{e}^{-( \tau-\xi) A} Bu\left( \xi \right)  \mathrm{d} \xi - \tau \varphi_1(-\tau A) Bu_0.
\end{equation}
Once more, employing \eqref{parabolicsmoothing1} and \eqref{eq15}, the bound for $\mathrm{e}^{-(n-1) \tau A} \delta_{1}$ is given by 
\begin{equation*}
	\begin{aligned}
		\left \| \mathrm{e}^{-(n-1) \tau A} \delta_{1}  \right \| 
		& \leq \left \| \mathrm{e}^{-(n-1) \tau A} A^{\gamma} \right \|   \int_0^{ \tau} \left \| \mathrm{e}^{-( \tau-\xi) A} \right \| \left \| A^{-\gamma} B \right \| \left \| u\left( \xi \right) \right \| \mathrm{d} \xi   \\
		& \qquad + \tau \left \| \mathrm{e}^{-(n-1) \tau A} A^{\gamma} \right \| \left \| \varphi_1(-\tau A) \right \| \left \| A^{-\gamma}  B \right \| \left \| u_0 \right \|  \\
		&
		\leq  C t_{n-1}^{-\gamma} \int_0^{ \tau}     \mathrm{d} \sigma  +   C t_{n-1}^{-\gamma} \tau   
		\leq  C t_{n-1}^{-\gamma} \tau  .  	
	\end{aligned}
\end{equation*}

Henceforth, the estimate will be almost the same as in the estimates of \eqref{eq11} and \eqref{eq13}.~The only thing that changes is the estimate of $\sum_{j=0}^{n-1} \mathrm{e}^{-j \tau A} \delta_{n-j}$ with the application of Lemma \ref{lem7} instead of Lemma \ref{boundeduprime}. Therefore, we will only state the results of these estimates here. The term for \( j = 0 \) can be bounded 
\begin{equation*}
	\begin{aligned}
		\left \| \delta_{n}  \right \| 
		 \leq  C t_{n-1}^{-1} \int_0^{ \tau} \int_0^\xi (\tau-\xi)^{-\gamma }    \mathrm{d} \sigma   \mathrm{d} \xi   \leq C t_{n-1}^{-1} \tau^{2-\gamma}  . 
	\end{aligned}
\end{equation*}

whereas the remaining sum is bounded by 
{\allowdisplaybreaks
	\begin{align*}
		\left \| \sum_{j=1}^{n-2} \mathrm{e}^{-j \tau A} \delta_{n-j}  \right \| 
	 \leq C  \sum_{j=1}^{n-2}  \int_0^{ \tau} \int_0^\xi t_j^{-\gamma} t_{n-j-1}^{-1}    \mathrm{d} \sigma   \mathrm{d} \xi   \leq C \tau^{2-\zeta}  \sum_{j=1}^{n-2} t_j^{-\gamma } t_{n-j-1}^{-1+\zeta} \leq C t_n^{-\gamma+\zeta} \tau^{1-\zeta}   . 
	\end{align*}	
}Here and throughout the paper, $\zeta >0$ is any fixed small number. We now present the convergence result for the exponential Euler method.
\begin{thm}\label{theo2}
	Let Assumptions \ref{ass1} and \ref{ass2} be fulfilled.~Then the numerical solution of the initial value problem \eqref{4121} obtained with the exponential Euler method \eqref{Euler} satisfies the error bound 
			$$
	\left\|u_n-u\left(t_n\right)\right\| \leq C t_n^{-\gamma^+} \tau^{1^-} 
	$$		
			uniformly in $0 \leq n \tau \leq T$. Assuming \( u_0 \in \mathcal{D}(A) \), we have
	$$
\left\|u_n-u\left(t_n\right)\right\| \leq C \tau 
$$
		uniformly in $0 \leq n \tau \leq T$.	The constant $C$ depends on $T$, but it is independent of $n$ and $\tau$. 
\end{thm}

\begin{proof}
	Taking the norm in \eqref{eq95} and applying the triangle inequality, the estimate of $e_n$ has the following form 
	\begin{equation*}
		\left \| e_n \right \|    \leq C \tau  \sum_{j=1}^{n-2} t_{n-j}^{-\gamma }   \left \| e_j \right \| +  C \tau^{1-\gamma} \left \| e_{n-1}     \right \|  + Ct_n^{-\gamma^+} \tau^{1^-},
\end{equation*}	
	wheares with assumption $u_0\in \mathcal{D}(A)$, we have
		\begin{equation*}
			\left \| e_n \right \|  \leq C \tau  \sum_{j=1}^{n-2} t_{n-j}^{-\gamma}   \left \| e_j \right \| +  C \tau^{1-\gamma} \left \| e_{n-1}     \right \|  + C\tau  .
	\end{equation*}
Applying a discrete Gronwall lemma (see \cite{HO2010}) concludes the proof.
\end{proof}
 
 \textit{Remark.} Theorem \ref{theo2} can also be obtained from \cite{HO2005} with a different derivation. With a convergence rate of $1^-$, there is no difference in the order of convergence from 1. The analysis in this section is intended solely for the motivating example, helping readers become familiar with the analysis when \( A \) and \( B \) are non-commuting and unbounded operators. The main results of this paper are presented in Sections \ref{secondorder} and \ref{thirdorder}. 

\subsection{Error bound for second-order methods}\label{secondorder}

This section addresses the convergence result for the unique one-parameter family of second-order exponential Runge--Kutta methods with $s=2$ in \eqref{eq20} (see \cite{HO2005}). 
\begin{remark}
	For the Taylor series expansion of \( u \), we distinguish between two cases: one where \( t_n = 0 \) and another where \( t_n \neq 0 \). The justification for this separation will be provided shortly. We get
	\begin{equation}\label{eq66} 
		u\left(\xi\right)=u\left(0\right)+\int_0^\xi u^{\prime}\left(\sigma\right) \mathrm{d} \sigma , \quad	u\left(t_n+\xi\right) = u\left(t_n\right)+ \xi u^{\prime} \left(t_n\right) + \int_0^\xi (\xi - \sigma) u^{\prime \prime}\left(t_n+\sigma\right) \mathrm{d} \sigma.
	\end{equation}
\end{remark}
\begin{lem}\label{lem4}
	Under Assumptions \ref{ass1}, \ref{ass2}, and supposing that $u_0\in \mathcal{D}(A)$, we have 
\begin{subequations}
	\begin{equation}\label{eq69}
		\left\| u^{\prime \prime}(t) \right\| \leq C t^{-1}.
	\end{equation} 					
	\begin{equation}\label{eq90}
		\left\| B u^{ \prime}(t) \right\| \leq C t^{-\gamma}.
	\end{equation}
\end{subequations}	
\end{lem} 
\begin{proof}
The bound of \eqref{eq69} follows in the same way as the estimate in Lemma \ref{boundeduprime}. The bound \eqref{eq90}, using the same technique as in \ref{boundeduprime}, can also be derived as follows
	\begin{equation}\label{eq70}
	\begin{aligned}
		\left\|  B u^{ \prime}(t)  \right\| 
		& \leq \left\| B A^{-\gamma} \right\|  \left\|A^{\gamma}  (A-B)^{-\gamma} \right\| \left\| (A-B)^{\gamma} \mathrm{e}^{-t(A-B)} \right\| \left\| (A-B)u_0 \right\|  \leq C t^{-\gamma}.
	\end{aligned}
\end{equation}
Here, the bound of $A^{\gamma}(A-B)^{-\gamma}$ can be obtained by observing that $A-B$ is only a slight perturbation of $A$.
\end{proof}
The integral part of \eqref{eq66} using Lemma \ref{lem4} can be bounded by
\begin{equation}\label{eq67}
	\left\| \int_0^\xi (\xi - \sigma) u^{\prime \prime}\left(t_n+\sigma\right) \mathrm{d} \sigma \right\| \leq C \int_0^\xi (\xi - \sigma) \left(t_n+\sigma\right)^{-1} \mathrm{d}\sigma \leq C \xi (\log(t_n+\xi)- \log(t_n)).	
\end{equation}
Here the value $t_n=0$ would make the integral in  \eqref{eq67} undefined. The separation \eqref{eq66} guarantees that two formulations in \eqref{eq66} are well-defined.

 From Section \ref{errorrecursion} with $s=2$, we have the following error expression for the unique one-parameter family of second-order exponential Runge--Kutta methods   \begin{equation}\label{eq17}
		e_n  = \text{I} + \text{II} + \sum_{l=1}^3 \text{III}_l, 
\end{equation}
where
\begin{equation*}
	\begin{aligned}
		\text{I} &= \tau \sum_{j=0}^{n-1} \mathrm{e}^{-(n-j-1) \tau A} \left( \sum_{l=1}^2 b_l(-\tau A) B \mathrm{e}^{-c_l \tau A} \right)  e_j, \\
		\text{II}  & = \tau^2 \sum_{j=0}^{n-1} \mathrm{e}^{-(n-j-1) \tau A} b_2\left(-\tau A\right) B   a_{21} \left(- \tau A\right )  B e_j, \, \\ \text{III}_l & = -\sum_{j=0}^{n-1} \mathrm{e}^{-j \tau A} \delta_{n-j}^{[l]}, \, 1\leq l \leq 3 .
	\end{aligned}
\end{equation*}		
		Here, we have
		\begin{equation*}
			\begin{aligned}
				\delta_{n}^{[1]} & = 
				 \tau b_2\left(-\tau A\right) B \int_0^{ c_2\tau} \mathrm{e}^{-\left( c_2\tau-\xi \right) A} \int_0^\xi Bu^{\prime}\left(t_{n-1}+\sigma\right) \mathrm{d} \sigma   \mathrm{d} \xi , \quad n\geq 1,\\
				\delta_{n}^{[2]} & = \int_0^{\tau} \mathrm{e}^{-( \tau-\xi) A} \int_0^\xi (\xi - \sigma) Bu^{\prime \prime}\left(t_{n-1}+\sigma\right) \mathrm{d} \sigma  \mathrm{d} \xi, \quad n\geq 2,\\
				\delta_{n}^{[3]} & = - \tau b_2\left(-\tau A\right)\int_0^{c_2\tau} \left(c_2\tau - \sigma \right) Bu^{\prime \prime}\left(t_{n-1}+\sigma\right) \mathrm{d} \sigma , \quad n\geq 2,
			\end{aligned}
		\end{equation*}  
		\begin{equation*}
			\begin{aligned}
				\delta_{1}^{[2]} & = \int_0^{\tau} \mathrm{e}^{-( \tau-\xi) A} \int_0^\xi Bu^{ \prime}\left(\sigma\right) \mathrm{d} \sigma   \mathrm{d} \xi, \qquad
				\delta_{1}^{[3]}  = -\tau b_2\left(-\tau A\right)\int_0^{c_2\tau} B u^{ \prime}\left(\sigma\right) \mathrm{d} \sigma.
			\end{aligned}
		\end{equation*}  
		As a result of \eqref{eq66}, we obtain two different expressions for $\delta_1^{[2]}$, $\delta_1^{[3]}$, $\delta_n^{[2]}$, and $\delta_n^{[3]}$, $n \geq 2$.
Henceforth, we assume that $u_0\in \mathcal{D}(A)$.~We start the estimation with a supporting result.
		\begin{lem}\label{lem3}
			Under Assumptions \ref{ass1} and \ref{ass2}, we have
			\begin{subequations}
				\begin{equation}\label{eq68a}
					\left\|   \phi \left(-\tau A\right) B\right\|\leq C\tau^{-\gamma},
				\end{equation} 					
				\begin{equation}\label{eq68b}
					\left\|  \mathrm{e}^{-t A} \phi \left(-\tau A\right) B\right\|\leq Ct^{-\gamma},
				\end{equation}

			\end{subequations}
			where $\phi(-\tau A)$ is an arbitrary linear combination of the functions $\varphi_k(-\tau A)$.
		\end{lem}
		\begin{proof}
We add the component \( A^{\gamma}A^{-\gamma} \) to each of the formulas \eqref{eq68a} and \eqref{eq68b}. The proof is completed using \eqref{parabolicsmoothing1} and \eqref{eq15}.
		\end{proof}
		The bound of I can be obtained using \eqref{eq68a}, \eqref{eq68b}, and \eqref{parabolicsmoothing1} as follows
		\begin{equation}\label{er6} 
			\begin{aligned}
				\left \| \text{I}  \right \| 
				& \leq    C \tau \sum_{j=0}^{n-2} \left \| \mathrm{e}^{-(n-j-1) \tau A}  \sum_{l=1}^2 b_l(-\tau A)  B   \right\| \left \| e_j \right \| +  C \tau \left \| \sum_{l=1}^2 b_l(-\tau A)  B \right\|   \left \| e_{n-1}      \right \| \\
				& \leq   C \tau \sum_{j=0}^{n-2} t_{n-j-1}^{-\gamma}   \left \| e_j \right \| +  C \tau^{1-\gamma} \left \| e_{n-1}     \right \|.					   
			\end{aligned}
		\end{equation}	
		The bound of II can be derived in a similar way with I, we omit the detail here. Lastly, we bound III.~The term for $j=0$ can be bounded using \eqref{parabolicsmoothing1}, \eqref{eq68a}, Lemmas \ref{boundeduprime}, and \ref{lem4}
		{\allowdisplaybreaks
			\begin{align*}
				\left \|	\delta_{n}^{[1]} \right \| 
& \leq 	 \tau	\left \|   b_2\left(-\tau A\right) B \right \| \int_0^{ c_2\tau} \int_0^\xi \left \| \mathrm{e}^{- \left( c_2\tau-\xi \right) A}  \right \|  \left \| B u^{\prime}\left(t_{n-1} +\sigma\right) \right \| \mathrm{d} \sigma   \mathrm{d} \xi    \\ 
& 				 \leq 	C t_{n-1}^{-\gamma} \sup_{0 \leq t \leq T} \left\|  u^{\prime}(t) \right\| \tau^{1-\gamma} 	\int_0^{ c_2 \tau}  \int_0^\xi   \mathrm{d} \sigma   \mathrm{d} \xi   \leq 	C t_{n-1}^{-\gamma} \tau^{3-\gamma}	 	.
			\end{align*}
		}Next, by employing \eqref{eq68a} and Lemma \ref{lem4}, we get
		\begin{equation*}\label{er5}
			\begin{aligned}
				\left \|	\delta_{n}^{[2]} \right \| 
				& \leq   \int_0^{\tau} \int_0^\xi \left \| \mathrm{e}^{-( \tau-\xi) A} B \right \| (\xi - \sigma)  \left \| u^{\prime \prime}\left(t_{n-1}+\sigma\right) \right \| \mathrm{d} \sigma   \mathrm{d} \xi   \\
				& \leq 	C t_{n-1}^{-1}	\int_0^{ \tau} \int_0^\xi \left( \tau-\xi \right)^{-\gamma} (\xi - \sigma)  \mathrm{d} \sigma   \mathrm{d} \xi   \leq 	C t_{n-1}^{-1} \tau^{3-\gamma}	.	
			\end{aligned}
		\end{equation*}
		The last term \( \delta_{n}^{[3]}, \, n \geq 0 \)  can be bounded in a similar way with \( \delta_{n}^{[2]},\, n \geq 0 \). The last term of the sum also needs to be treated separately. The terms for $j=n-1$ can be bounded by utilising \eqref{parabolicsmoothing1}, \eqref{eq15}, \eqref{eq68a}, \eqref{eq68b}, and Lemma \ref{boundeduprime} as follows
		\begin{equation*}\label{er5}
			\begin{aligned}
				\left \| \mathrm{e}^{-(n-1) \tau A}	\delta_{1}^{[1]} \right \| 
				& \leq  \tau \left \| \mathrm{e}^{-(n-1) \tau A}  b_2\left(-\tau A\right)  B \right \| \int_0^{ c_2\tau} \int_0^\xi \left \| \mathrm{e}^{-\left( c_2\tau-\xi \right) A} B \right \|  \left \| u^{\prime}\left(\sigma\right) \right \| \mathrm{d} \sigma   \mathrm{d} \xi  \\ & \leq 	C t_{n-1}^{-\gamma} \sup_{0 \leq t \leq T} \left\|  u^{\prime}(t) \right\| \tau  	\int_0^{ c_2\tau} \int_0^\xi  \left( c_2\tau-\xi \right)^{-\gamma}  \mathrm{d} \sigma   \mathrm{d} \xi   \leq 	C t_{n-1}^{-\gamma}  \tau^{3-\gamma}	 .	
			\end{aligned}
		\end{equation*}
		For the second one, we get
		\begin{equation*}\label{er5}
			\begin{aligned}
				\left \|\mathrm{e}^{-(n-1) \tau A}	\delta_{1}^{[2]} \right \|
				& \leq  \left \| \mathrm{e}^{-(n-1) \tau A} A^{\gamma}   \right \|  \int_0^{\tau} \int_0^\xi \left \| \mathrm{e}^{-( \tau-\xi) A} \right \|    \left \| A^{-\gamma}B  \right \|    \left \| u^{ \prime}\left(\sigma\right) \right \| \mathrm{d} \sigma   \mathrm{d} \xi  \\		
				 & \leq 	C t_{n-1}^{-\gamma} \sup_{0 \leq t \leq T} \left\|  u^{\prime}(t) \right\|	\int_0^{ \tau} \int_0^\xi     \mathrm{d} \sigma   \mathrm{d} \xi  \leq 	C t_{n-1}^{-\gamma} \tau^{2},
			\end{aligned}
		\end{equation*}

		The remaining sums with $j\neq 0$ and $j\neq n-1$ can be bounded by using \eqref{parabolicsmoothing1}, \eqref{eq15}, \eqref{eq68a}, \eqref{eq68b}, and Lemma  \ref{lem4}  as follows
		\begin{equation*}\label{er5}
			\begin{aligned}
				\left \|	\sum_{j=1}^{n-2} \mathrm{e}^{-j h A} \delta_{n-j}^{[1]} \right \| 
				\leq &	\tau	 \sum_{j=1}^{n-2} \left \| \mathrm{e}^{-j \tau A}     b_2\left(-\tau A\right) B \right \| \int_0^{ c_2\tau} \int_0^\xi \left \| \mathrm{e}^{-\left( c_2\tau-\xi \right) A}  \right \|  \left \| Bu^{\prime}\left(t_{n-j-1} +\sigma\right) \right \| \mathrm{d} \sigma   \mathrm{d} \xi  \\ 				  \leq  &	C\tau  \sum_{j=1}^{n-2} t_{j}^{-\gamma} t_{n-j-1}^{-\gamma}  \int_0^{ c_2\tau} \int_0^\xi    \mathrm{d} \sigma   \mathrm{d} \xi  
				\leq 	C t_n^{-2\gamma+1} \tau^{2},
			\end{aligned}
		\end{equation*}
		\begin{equation*}\label{er5}
			\begin{aligned}
				\left \|	\sum_{j=1}^{n-2} \mathrm{e}^{-j \tau A} \delta_{n-j}^{[2]} \right \| 
				& \leq  C \sum_{j=1}^{n-2}\left \| \mathrm{e}^{-j \tau A} A^{\gamma} \right \| \int_0^{\tau} \int_0^\xi  (\xi - \sigma) \left \| \mathrm{e}^{-\left( \tau-\xi \right) A} \right \|  \left \| A^{-\gamma}B \right \|  \left \| u^{\prime \prime}\left(t_{n-j-1}+\sigma\right) \right \|  \mathrm{d} \sigma   \mathrm{d} \xi  \\		&				 \leq   C \sum_{j=1}^{n-2}	t_{j}^{-\gamma} t_{n-j-1}^{-1} \int_0^{ \tau} \int_0^\xi  (\xi - \sigma)  \mathrm{d} \sigma   \mathrm{d} \xi  
				\leq 	C t_n^{-\gamma+\zeta} \tau^{2-\zeta}. 
			\end{aligned}
		\end{equation*} 	
		  Now, we are ready to present the convergence result.
		\begin{thm}\label{theo3}
			Let Assumptions \ref{ass1} and \ref{ass2} be fulfilled. In addition, we assume that $u_0\in \mathcal{D}(A)$.~Then the numerical solution of the initial value problem \eqref{4121} using exponential Runge--Kutta methods of second-order with $s=2$ in \eqref{eq20} satisfies the error bound
			$$
			\left\|u_n-u\left(t_n\right)\right\| \leq C t_n^{-\gamma^+} \tau^{2^-} 
			$$
			uniformly in $0 \leq n \tau \leq T$. The constant $C$ depends on $T$, but it is independent of $n$ and $\tau$.
		\end{thm}
		
		\begin{proof}
			By taking the norm in \eqref{eq17} and utilizing the triangle inequality, the estimate of $e_n$ has the following form
			\begin{equation*}
					\left \| e_n \right \|   \leq C \tau  \sum_{j=1}^{n-2} t_{n-j}^{-\gamma}   
					\left \| e_j \right \|  + C \tau^{1-\gamma}    \left \| e_{n-1} \right \|      + C t_n^{-\gamma^+} \tau^{2^-}. 
			\end{equation*}
			An application of a discrete Gronwall lemma (see \cite{HO2010}) concludes the proof. 
		\end{proof}

		\subsection{Error bound for third-order methods}\label{thirdorder} 
		
		We now move on to the main part of this paper, the error analysis of third-order exponential Runge--Kutta methods with $s=3$ in \eqref{eq20} (see \cite{HO2005}).~We assume that there exist some \( 0 < \alpha \leq \frac{1}{2} \) for which \(u_0 \in \mathcal{D}(A^{1+\alpha}) \) in this section (see \cite{HOCHBRUCK2005323}). 

		 \begin{remark}
		 	Again, we distinguish two cases: $t_n = 0$ and $t_n\neq 0$ as in \eqref{eq66}.~By expanding $u$ in a Taylor series, we obtain
		 	\begin{equation}\label{errf32}
		 		u\left( \xi\right) = u\left( 0 \right)+ \xi u^{\prime} \left( 0 \right)  + \int_0^\xi (\xi - \sigma) u^{ \prime \prime}\left(\sigma\right) \mathrm{d} \sigma ,
		 	\end{equation}
		 	and 
		 	\begin{equation}\label{eq59} 
		 		u\left( t_n+\xi\right) = u\left( t_n \right)+ \xi u^{\prime} \left( t_n \right)  + \frac{\xi^2}{2} u^{ \prime \prime}(t_n) + \int_0^\xi \frac{(\xi - \sigma)^2}{2} u^{ \prime \prime \prime}\left(t_n+\sigma\right) \mathrm{d} \sigma .
		 	\end{equation}
		 	The expression for the second- and third-order derivative of 
		 	$u$ can be formulated as 
		 	\begin{equation}\label{eq16}
		 		u^{\prime \prime }(\sigma) =(A-B)^{1-\alpha} \mathrm{e}^{-\sigma(A-B)} (A-B)^{1+\alpha}u_0,\, \, \, \,  u^{\prime \prime \prime}(t) =(A-B)^{2-\alpha} \mathrm{e}^{-t(A-B)} (A-B)^{1+\alpha}u_0.
		 	\end{equation}

		 	\begin{lem}\label{lem5} 
		 		Under Assumptions \ref{ass1}, \ref{ass2},  and \( u_0 \in \mathcal{D}(A^{1+\alpha}) \), we have 
\begin{subequations}
	\begin{equation}\label{eq73}
		 			\left\| u^{\prime \prime }(t) \right\|  \leq C t^{-1+\alpha}, 
	\end{equation} 					
	\begin{equation}\label{eq74}
		 			\left\| u^{\prime \prime \prime}(t) \right\| \leq C t^{-2+\alpha}. 
	\end{equation}
\end{subequations}			 		
		 	\end{lem} 
		 	\begin{proof}
		 		These bounds follow at once by employing a property similar to  \eqref{parabolicsmoothing1} and \( u_0 \in \mathcal{D}(A^{1+\alpha}) \).
		 	\end{proof}
		 	
		 	Applying Lemma \ref{lem5}, the integral expressions in  \eqref{eq59} can be bounded by
		 	
		 	\begin{equation}\label{eq61}
		 		\left\| \int_0^\xi \frac{(\xi - \sigma)^2}{2} u^{\prime \prime \prime}\left(t_n+\sigma\right) \mathrm{d} \sigma \right\| \leq C \int_0^\xi \frac{(\xi - \sigma)^2}{2} \left(t_n+\sigma\right)^{-2+\alpha} \mathrm{d} \sigma \leq C t_n^{-1+\alpha} \xi^2 .	
		 	\end{equation}
		 	Here \eqref{eq61} would be undefined with the value $t_n=0$, which explains \eqref{errf32}. 
		 \end{remark}
		
		Assuming the order conditions No. 1-4 in \eqref{Table1} are satisfied in a strong form, whereas No. 5 in \eqref{Table1} is fulfilled in a very weak form with $A=0$ in Table \ref{Table1}. Taking $s=3$ and following the procedure in Section \ref{errorrecursion}, we obtain the expression for the global error 
	\begin{equation}\label{eq18}
		e_n
		 = \text{I} + \text{II} + \text{III}+ \text{IV} + \text{V} + \sum_{l=1}^7 \text{VI}_l,
\end{equation}
where
{\allowdisplaybreaks
		\begin{align}
				\text{I} & =  \tau \sum_{j=0}^{n-1} \mathrm{e}^{-(n-j-1) \tau A} \left( \sum_{l=1}^3 b_l(-\tau A) B \mathrm{e}^{-c_l \tau A} \right) e_j, \, \nonumber \\
				\text{II} & =  \tau^2 \sum_{j=0}^{n-1} \mathrm{e}^{-(n-j-1) \tau A} 		 \left( b_2(-\tau A) B      a_{21}\left(-\tau A\right)  B + \sum_{l=1}^2 b_3(-\tau A) B a_{3l}\left(-\tau A\right) B \mathrm{e}^{-c_l \tau A} \right) e_j, \nonumber \\
				\text{III} & =  \tau^3 \sum_{j=0}^{n-1} \mathrm{e}^{-(n-j-1) \tau A}   b_3(-\tau A) B      a_{32}\left( -\tau A\right) B a_{21} \left(- \tau A\right) B e_j, \, \nonumber \\ 				
				\text{IV} & = - \tau^3 \sum_{j=0}^{n-1} \mathrm{e}^{-j \tau A} \sum_{l=1}^2   \mathcal{T}_{n-j}^{[l]}, \nonumber \\  
				\text{V} & = -\tau^4 \sum_{j=0}^{n-1} \mathrm{e}^{-j \tau A}  \mathcal{L}_{n-j} , \nonumber \\				   
			    \text{VI}_l & = - \sum_{j=0}^{n-1} \mathrm{e}^{-j \tau A}  \mathcal{R}_{n-j}^{[l]},  1\leq l\leq 7. \nonumber 
			\end{align}
		}
		Here, we have
		{\allowdisplaybreaks
			\begin{align*}
				\mathcal{T}_{n}^{[1]} & =   c_2 b_2(-\tau A) B    a_{21} \left(- \tau A\right) Bu^{\prime}\left(t_{n-1}\right) , \qquad n \geq 1,
				\\
				\mathcal{T}_{n}^{[2]} & = b_3(-\tau A) B      \left( c_3^2 \varphi_2(-c_3 \tau A)-c_2 a_{32}(-\tau A) \right) Bu^{\prime}\left(t_{n-1}\right) , \qquad n \geq 1,\\
				\mathcal{L}_{n} & =  - c_2 b_3(-\tau A) B    a_{32}\left( -\tau A\right)B a_{21} \left(- \tau A\right) Bu^{\prime}\left(t_{n-1}\right), \qquad n \geq 1, \\				
				\mathcal{R}_{n}^{[k]}	& =   \tau b_{k+1}(-\tau A)B   \int_0^{c_{k+1} \tau} \mathrm{e}^{-(c_{k+1} \tau-\xi) A} \int_0^\xi (\xi - \sigma) Bu^{\prime \prime}\left(t_{n-1}+\sigma\right) \mathrm{d} \sigma   \mathrm{d} \xi , \quad n \geq 1, \quad  k =1,2,
				\\
				\mathcal{R}_{n}^{[3]} & =  \tau^2 b_3(-\tau A) B    a_{32}\left( -\tau A\right) B   \int_0^{c_2 \tau} \mathrm{e}^{-(c_2 \tau-\xi) A} \int_0^\xi (\xi - \sigma) Bu^{\prime \prime}\left(t_{n-1}+\sigma\right) \mathrm{d} \sigma   \mathrm{d} \xi ,  \qquad n \geq 1,  \\
				\mathcal{R}_{n}^{[4]} & = - \tau^2 b_3(-\tau A) B  a_{32}\left(-\tau A\right) B \int_0^{c_2\tau} \left(c_2\tau - \sigma\right) Bu^{\prime \prime}\left(t_{n-1}+\sigma\right) \mathrm{d} \sigma  ,   \qquad n \geq 1,  \\
				\mathcal{R}_{n}^{[5]} & = \int_0^{ \tau} \mathrm{e}^{-(\tau-\xi) A} \int_0^\xi \frac{ (\xi - \sigma)^2}{2}  Bu^{\prime \prime \prime}\left(t_{n-1}+\sigma\right) \mathrm{d} \sigma   \mathrm{d} \xi, \qquad n \geq 2, \\
				\mathcal{R}_{n}^{[m]} & = - \tau b_{m-4}\left(-\tau A\right) \int_0^{c_{m-4}\tau} \frac{\left(c_{m-4}\tau - \sigma \right)^2}{2} Bu^{\prime \prime \prime}\left(t_{n-1}+\sigma\right) \mathrm{d} \sigma , \qquad n \geq 2, \quad  m =6,7, \\
				\mathcal{R}_{1}^{[5]} & = \int_0^{ \tau} \mathrm{e}^{-(\tau-\xi) A} \int_0^\xi  (\xi - \sigma)  Bu^{\prime  \prime}\left(\sigma\right) \mathrm{d} \sigma   \mathrm{d} \xi, \\ 
\mathcal{R}_{1}^{[m-4]} & = - \tau b_{m-4}\left(-\tau A\right) \int_0^{c_{m-4}\tau} (c_{m-4}\tau - \sigma ) Bu^{ \prime \prime}\left(\sigma\right) \mathrm{d} \sigma,		\quad  m =6,7.			
			\end{align*}
		}We start deriving the error bound for  $e_n$. We observe that the estimates for the terms in parentheses in both $\text{I}$ and $\text{II}$ are similar. Furthermore, the estimation of I, II, and III can be done in the same way as in \eqref{er6}. Therefore, we only show the result here
		{\allowdisplaybreaks
			\begin{align*}
				\left \| \text{I}  \right \| 
				& \leq   C \tau \sum_{j=0}^{n-2} t_{n-j-1}^{-\gamma}   \left \| e_j \right \| +   C \tau^{1-\gamma} \left \| e_{n-1}     \right \| , \quad
				\left \| \text{II}  \right \| 
				 \leq C \tau \sum_{j=0}^{n-1}\left\|e_j\right\|, \quad
				\left \| \text{III}  \right \| 
				 \leq C \tau \sum_{j=0}^{n-1}\left\|e_j\right\|.     
			\end{align*}
		}The following identities turn out to be crucial  for estimating $ \text{IV} $. There is a bounded operator denoted as \(\widetilde{\phi}(-\tau A)\) with the following property (see \cite{HO2005})
		\begin{align}\label{eq55}
			\phi(-\tau A)-\phi(0) & = (\tau A)^{\widehat{\Gamma}} \widetilde{\phi}(-\tau A),   \qquad 0 \leq \widehat{\Gamma} < 1.
		\end{align}
		\begin{lem}\label{lem8}
			Under Assumptions \ref{ass1}, \ref{ass2}, $\left \| A^{\frac{1}{2}\widetilde{\Gamma}} BA^{-1} \right \|\leq C $ where $\widetilde{\Gamma} \in \mathbb{R}$ is fixed, and $u_0 \in \mathcal{D}(A^{1+\alpha})$, the following bound holds
			\begin{equation}\label{eq77}
				\left\|  B \phi\left(-\tau A\right) B u^{\prime}\left(t\right) \right\|\leq \tau^{\frac{1}{2}\widetilde{\Gamma}-\gamma}t^{-1+\alpha},
			\end{equation} 					
			where $\phi(-\tau A)$ is an arbitrary linear combination of the functions $\varphi_k(-\tau A)$.
		\end{lem}
		\begin{proof}
			We can obtain the bound as follows, 
			\begin{align*}
				\left\|  B \phi\left(-\tau A\right) B u^{\prime}\left(t\right) \right\| 
				& \leq \left \|  B A^{-\gamma} \right \| \left \|\phi\left(-\tau A\right)A^{-\frac{1}{2}\widetilde{\Gamma}+\gamma} \right \| \left \| A^{\frac{1}{2}\widetilde{\Gamma}} BA^{-1} \right \| \left \|A (A-B)^{-1} \right \| \left \| u^{\prime \prime}\left(t\right) \right \| \\
				& \leq \tau^{\frac{1}{2}\widetilde{\Gamma}-\gamma}t^{-1+\alpha}.
			\end{align*}
Here, the explanation is similar to that in Lemma \ref{lem4}.
		\end{proof}
		 Employing Lemma \ref{lem8}, the term for $j = 0$ can be bounded by
		{\allowdisplaybreaks
			\begin{align*}
				\left \| \mathcal{T}_{n}^{[1]} + \mathcal{T}_{n}^{[2]} \right \|  
				& \leq  C \tau^3  \left \|   b_2(-\tau A) \right \| \left \| B  a_{21} \left(- \tau A\right) B  u^{\prime}\left(t_{n-1}\right) \right \| \\
				& \qquad + C \tau^3  \left \|    b_3(-\tau A) \right \| \left \| B      \left( c_3^2 \varphi_2(- \tau A)-c_2 a_{32}(-\tau A) \right)    B  u^{\prime}\left(t_{n-1}\right) \right \| \\
				&  \leq C t_{n-1}^{-1+\alpha} \tau^{3 + \frac{1}{2}\widetilde{\Gamma}-\gamma} .
			\end{align*}
		}In the same way, the term for $j = n-1$ can be bounded by	
		{\allowdisplaybreaks
			\begin{align*}
				 \left \| \mathrm{e}^{-(n-1) \tau A}  \mathcal{T}_{1}^{[1]} + \mathrm{e}^{-(n-1) \tau A}  \mathcal{T}_{1}^{[2]} \right \|  
				& \leq  C \tau^3   \left \| A^{-1} B A^{\frac{1}{2}\widetilde{\Gamma}} \right \| \left \|  A^{-\frac{1}{2}\widetilde{\Gamma} +\frac{1}{2}}  a_{21} \left(- \tau A\right) \right \| \\
				& \quad + C \tau^3   \left \| A^{-1}   B   A^{\frac{1}{2}\widetilde{\Gamma}} \right \|  \left \| A^{-\frac{1}{2}\widetilde{\Gamma} +\frac{1}{2}} \left( c_3^2 \varphi_2(- \tau A)-c_2 a_{32}(-\tau A) \right)  \right \|   \\			
				&  \leq C t_{n-1}^{-1+\alpha} \tau^{3 + \frac{1}{2}\widetilde{\Gamma}-\gamma} .		
			\end{align*}
		}Here, we additionally assume that \(\left \| A^{-1}  B A^{\frac{1}{2}\widetilde{\Gamma}} \right \| \leq C\) and $\left \| A^{\frac{1}{2}\widetilde{\Gamma}} B A^{-1}  \right \| \leq C $, where \(\widetilde{\Gamma} \in \mathbb{R}  \) is fixed. These assumptions will be discussed in more detail later.
	
Using condition No. 5 with all arguments evaluated for $A=0$, and employing \eqref{eq55} with $\phi = b_2(-\tau A), b_3(-\tau A), \varphi_2(-\tau A), \varphi_3(-\tau A)$, the sum of the remaining terms with $j \neq 0$ and $j \neq n-1$ can be expressed as follows
		\begin{equation}\label{eq21}
						\begin{aligned}
				&  \tau^3 \sum_{j=1}^{n-2} \mathrm{e}^{-j \tau A}   \mathcal{T}_{n-j}^{[1]}  + \tau^3 \sum_{j=1}^{n-2} \mathrm{e}^{-j \tau A} \mathcal{T}_{n-j}^{[2]}  
				\\
				& =  \tau^{3+\Gamma} \sum_{j=1}^{n-2}  \left( \mathrm{e}^{-j \tau A}   A^{\Gamma } \right)  \widetilde{b_2}(-\tau A)   \left ( B  c_2 a_{21} \left(- \tau A\right) B u^{\prime }\left(t_{n-j-1}\right) \right )  \\
				& \quad	+	 \tau^{3+\Gamma} \sum_{j=1}^{n-2} \left ( \mathrm{e}^{-j \tau A}    A^{\Gamma } \right )  \widetilde{b_3}(-\tau A)  \left ( B \left(   c_3^2 \varphi_2(- \tau A)-c_2 a_{32}(-\tau A)   \right)   B u^{\prime }\left(t_{n-j-1}\right) \right )  \\
				& \quad +  \tau^{3 + \widetilde{\Gamma}}  \sum_{j=1}^{n-2} \left ( \mathrm{e}^{-j \tau A} A \right )   A^{-1} B  A^{\frac{1}{2}\widetilde{\Gamma}}  \left ( c_2 \widetilde{\varphi_2}(-\tau A) \right )  A^{\frac{1}{2}\widetilde{\Gamma}} B A^{-1}  \left ( A (A-B)^{-1} \right )  u^{\prime \prime }\left(t_{n-j-1}\right)  \\
				& \quad	+  \tau^{3 + \widetilde{\Gamma}} \sum_{j=1}^{n-2} \left ( \mathrm{e}^{-j \tau A}   A  \right )  A^{-1}  B A^{\frac{1}{2}\widetilde{\Gamma}} \left (  \widetilde{\varphi_3}(-\tau A) \right )  A^{\frac{1}{2}\widetilde{\Gamma}} B A^{-1}  \left ( A (A-B)^{-1} \right )  u^{\prime \prime}\left(t_{n-j-1}\right) . 	\end{aligned} 	
		\end{equation}
		Using \eqref{parabolicsmoothing1}, and Lemmas \ref{lem4}, \ref{lem8} to estimate each term of the parentheses in \eqref{eq21}, we obtain
		{\allowdisplaybreaks
	\begin{align*}
		  \left \|  \tau^3 \sum_{j=1}^{n-2} \mathrm{e}^{-j \tau A}   \mathcal{T}_{n-j}^{[1]}  + \tau^3 \sum_{j=1}^{n-2} \mathrm{e}^{-j \tau A} \mathcal{T}_{n-j}^{[2]}  \right \|
		& \leq C \tau^{ 3 + \Gamma+ \frac{1}{2}\widetilde{\Gamma}-\gamma}  \sum_{j=1}^{n-2}	t_{j}^{- \Gamma } t_{n-j-1}^{-1+\alpha} + C \tau^{3 + \widetilde{\Gamma}}  \sum_{j=1}^{n-2}	t_{j}^{-1} t_{n-j-1}^{-1+\alpha} \\
		& \leq C t_n^{-\Gamma + \alpha} \tau^{2 + \Gamma+ \frac{1}{2}\widetilde{\Gamma}-\gamma } + C t_n^{-1 + \alpha +\zeta } \tau^{2 + \widetilde{\Gamma}-\zeta}.	\end{align*}		
}By taking $\Gamma = 1- \zeta,\text{ } \widehat{\Gamma} =  \widetilde{\Gamma} =   \frac{1}{2}- \zeta$  in \eqref{eq55}, we obtain the order of convergence as $ \frac{5}{2}-2\zeta$ for the estimation of IV.
		
We continue to bound V, with the application of \eqref{eq68a}, \eqref{eq68b}, and Lemmas \ref{boundeduprime}, \ref{lem8}. The term for \( j = 0 \) can be bounded by
		\begin{equation*}
			\begin{aligned}
				\left \| \tau^4 \mathcal{L}_{n} \right \| 
				\leq  C \tau^4  \left \|   b_3(-\tau A)  B \right \| \left \|    a_{32}\left( -\tau A\right)  \right \| \left \| B a_{21} \left(- \tau A\right)  Bu^{\prime}\left(t_{n-1}\right)    \right \|  \leq C t_{n-1}^{-1+\alpha} \tau^{4 + \frac{1}{2}\widetilde{\Gamma}- 2\gamma}.
			\end{aligned}
		\end{equation*}
		The term for $j = n-1$ can be estimated by	
		\begin{equation*}
			\begin{aligned}
				\left \| \tau^4 \mathrm{e}^{-(n-1) \tau A} \mathcal{L}_{1} \right \| 
				& \leq  C \tau^4  \left \|  \mathrm{e}^{-(n-1) \tau A} b_3(-\tau A)  B \right \| \left \|    a_{32}\left( -\tau A\right) B  \right \| \left \|  a_{21} \left(- \tau A\right)  B \right \| \left \| u^{\prime}\left(0\right)    \right \| \\ & \leq C t_{n-1}^{-\gamma} \tau^{4-2\gamma},
			\end{aligned}
		\end{equation*}
		whereas the sum of remaining terms with $j \neq 0$ and $j \neq n-1$ can be bounded by
		\begin{equation*}
			\begin{aligned}
				\left \| \tau^4 \sum_{j=1}^{n-2} \mathrm{e}^{-j \tau A}  \mathcal{L}_{n-j} \right \| 
								& \leq  C \tau^4 \sum_{j=1}^{n-2} \left \| \mathrm{e}^{-j \tau A}   b_3\left(-\tau A\right)     B \right \| \left \|     a_{32}\left( -\tau A\right)  \right \| \left \| B a_{21} \left(- \tau A\right)  Bu^{\prime}\left(t_{n-j-1}\right)      \right \| \\ 
				 & \leq C \tau^{4 + \frac{1}{2}\widetilde{\Gamma}-\gamma}  \sum_{j=1}^{n-2}	t_{j}^{ -\gamma}t_{n-j-1}^{-1+\alpha}   \leq C t_n^{-\gamma+\alpha} \tau^{3+ \frac{1}{2}\widetilde{\Gamma}-\gamma }  . 		
			\end{aligned}
		\end{equation*}
		By the same choice of $\widehat{\Gamma} =  \widetilde{\Gamma} =   \frac{1}{2}- \zeta$ in \eqref{eq55}, we get the convergence rate of $\frac{13}{4} - \frac{1}{2} \zeta -\gamma $ for V.
		
		Finally, we bound  $\text{VI}$.~We notice that the estimations of $\mathcal{R}_{n}^{[1]}$ and $\mathcal{R}_{n}^{[2]}$, as well as $\mathcal{R}_{n}^{[6]}$ and $\mathcal{R}_{n}^{[7]}$, respectively, are the same.~Therefore, we estimate only one.~Using \eqref{eq68a} and \eqref{eq73}, the first two terms for $j=0$ can be bounded by
		{\allowdisplaybreaks
			\begin{align*}
				\left \|   \mathcal{R}_{n}^{[1]}  \right \| 
				& \leq C  \tau \left \| b_2(-\tau A) B  \right \|  \int_0^{c_2 \tau} \int_0^\xi \left \| \mathrm{e}^{-\left( c_2 \tau-\xi \right) A} B    \right \| (\xi - \sigma)  \left \| u^{ \prime \prime}\left(t_{n-1}+\sigma\right) \right \| \mathrm{d} \sigma   \mathrm{d} \xi  \\								
				& \leq 	C t_{n-1}^{\alpha-1} \tau^{1 - \gamma} 	\int_0^{c_2 \tau} \int_0^\xi  \left( c_2 \tau-\xi \right)^{- \gamma} (\xi - \sigma)  \mathrm{d} \sigma   \mathrm{d} \xi  \leq 	C t_{n-1}^{\alpha-1} \tau^{4- 2\gamma} ,\\
				\left \| \mathcal{R}_{n}^{[3]} \right \|
				& \leq C \tau^2 \left\|  b_3(-\tau A)  B  \right \| \left\|  a_{32}\left( -\tau A\right) B \right \| \times \\ & \qquad \int_0^{c_2 \tau} \int_0^\xi  \left \| \mathrm{e}^{-\left( c_2 \tau-\xi \right) A} B    \right \|   (\xi - \sigma)  \left \| u^{ \prime \prime}\left(t_{n-1}+\sigma\right)  \right \|  \mathrm{d} \sigma   \mathrm{d} \xi  \\				  
				& \leq 	C t_{n-1}^{\alpha - 1} \tau^{2-2\gamma} 	\int_0^{c_2 \tau} \int_0^\xi \left( c_2 \tau-\xi \right)^{-\gamma} (\xi - \sigma)  \mathrm{d} \sigma   \mathrm{d} \xi  \leq 	C t_{n-1}^{\alpha - 1} \tau^{5-3\gamma}	.	 		
			\end{align*}
		}We observe that the estimate of $\mathcal{R}_{n}^{[3]}$ is almost the same as $\mathcal{R}_{n}^{[1]}$. The only difference is that the order of $\mathcal{R}{n}^{[3]}$ will be $1 - \gamma$ order higher than that of $\mathcal{R}{n}^{[1]}$ due to the additional term $\tau \left \| a_{32}\left( -\tau A\right) B \right \| \leq \tau^{1-\gamma}$.
		For this reason, we will only estimate related to $\mathcal{R}_{l}^{[1]}$, $1 \leq l \leq n-1 $.~In the same way, we proceed further with the other terms
		\begin{equation*}
			\begin{aligned}
				\left \| \mathcal{R}_{n}^{[4]}  \right \| 
				& \leq \tau^2    \left \| b_3(-\tau A) B \right \| \left\|  a_{32}\left(-\tau A\right) B \right \| \int_0^{c_2\tau} \left(c_2\tau - \sigma \right)  \left\| u^{\prime \prime}\left(t_{n-1}+\sigma\right)  \right \| \mathrm{d} \sigma     \\				
				& \leq 	C t_{n-1}^{\alpha-1} \tau^{2-2\gamma} 	\int_0^{c_2\tau} \left(c_2\tau - \sigma \right)  \mathrm{d} \sigma  \leq 	C t_{n-1}^{\alpha-1} \tau^{4-2\gamma} 	.	
			\end{aligned}
		\end{equation*}
		Here, we notice that the result of the estimation of $\mathcal{R}_{n}^{[1]}$ and $\mathcal{R}_{n}^{[4]}$ is the same, so we will estimate them only once. This comes from the fact that the only main difference is that we have an additional integral in $\mathcal{R}_{n}^{[1]}$, whereas we have an additional $\tau$ in $\mathcal{R}_{n}^{[4]}$. This observation also holds for the estimation of $\mathcal{R}_{n}^{[5]}$ and $\mathcal{R}_{n}^{[6]}$, $\mathcal{R}_{n}^{[7]}$. The estimation of $\mathcal{R}_{n}^{[5]}$  use \eqref{eq68a} and \eqref{eq74} as follows
		{\allowdisplaybreaks
			\begin{align*}
				\left \| \mathcal{R}_{n}^{[5]}  \right \|
				& \leq    \int_0^{ \tau} \int_0^\xi   \left \| \mathrm{e}^{-\left(  \tau-\xi \right) A} B   \right \|  \frac{ (\xi - \sigma)^2}{2}  \left \| u^{\prime \prime \prime}\left(t_{n-1}+\sigma\right) \right \|  \mathrm{d} \sigma   \mathrm{d} \xi  \\				 
				& \leq 	C  t_{n-1}^{\alpha-2}	\int_0^{ \tau} \int_0^\xi  \left(  \tau-\xi \right)^{-\gamma}  \frac{ (\xi - \sigma)^2}{2}  \mathrm{d} \sigma   \mathrm{d} \xi  \leq 	C t_{n-1}^{\alpha-2} \tau^{4-\gamma} 	.	
			\end{align*}
		}The term for $j=n-1$ also needs to be estimated separately.~For the first one, we have
		\begin{equation*}
			\begin{aligned}
				\left \|\mathrm{e}^{-(n-1) \tau A}   \mathcal{R}_{1}^{[1]}  \right \| 
				& \leq \tau \left \| \mathrm{e}^{-(n-1) \tau A}  b_2(-\tau A)  B \right \|     \int_0^{c_2 \tau} \int_0^\xi  \left \|  \mathrm{e}^{-\left(c_2 \tau-\xi \right) A} B \right\|  (\xi - \sigma)   \left \| u^{ \prime \prime}\left(\sigma\right) \right \| \mathrm{d} \sigma   \mathrm{d} \xi  \\									
				& \leq 	C t_{n-1}^{ - \gamma} \tau 	\int_0^{c_2 \tau} \int_0^\xi \left( c_2 \tau-\xi \right)^{-\gamma}(\xi - \sigma) \sigma^{-1+\alpha}    \mathrm{d} \sigma   \mathrm{d} \xi   \leq 	C  t_{n-1}^{ - \gamma} \tau^{3+\alpha-\gamma} .		
			\end{aligned}
		\end{equation*}
			Here, we used \eqref{eq68a}, \eqref{eq68b}, and \eqref{eq73}. We establish the following lemma for estimating the remainders.
			\begin{lem}\label{lem6} 
				Under Assumptions \ref{ass1}, \ref{ass2}, and let $\beta \in \mathbb{R}$, we have   
\begin{subequations}
	\begin{equation}\label{eq71}
					\left\| (A-B)^{-\beta} u^{\prime \prime }(t) \right\|  \leq C t^{-1+\alpha+\beta}, 
	\end{equation} 					
	\begin{equation}\label{eq72}
					\left\| (A-B)^{-\beta} u^{\prime \prime \prime}(t) \right\| \leq C t^{-2+\alpha+\beta}. 
	\end{equation}
\end{subequations}			 						
			\end{lem} 
			\begin{proof}
				These bounds follow in the same manner as the estimate in Lemma \ref{lem5}.
			\end{proof}
			For the moment, we assume that $\left \| A^{-1}BA^\beta \right \|$ is bounded, which will be verified later.~At $t_n=0$, we use the Taylor expansion of 
			$u$  only up to the second-order \eqref{errf32} 
			\begin{equation*}
				\begin{aligned}
					&\left \| \mathrm{e}^{-(n-1) \tau A} \mathcal{R}_{1}^{[5]}  \right \|\leq  \left \| \mathrm{e}^{-(n-1) \tau A} A \right \| \times \\ & \qquad \int_0^{ \tau} \int_0^\xi \left\| \mathrm{e}^{-\left(  \tau-\xi \right) A} \right\|    (\xi - \sigma) \left \| A^{-1}BA^\beta \right \| \left \| A^{-\beta} (A-B)^{\beta}  \right \| \left \| (A-B)^{-\beta} u^{ \prime \prime}\left(\sigma\right) \right \|  \mathrm{d} \sigma   \mathrm{d} \xi  \\	&
					\qquad \qquad \qquad \quad \leq 	C  t_{n-1}^{ -1}	\int_0^{ \tau} \int_0^\xi  (\xi - \sigma) \sigma^{-1+\alpha+\beta}   \mathrm{d} \sigma   \mathrm{d} \xi  \leq 	C t_{n-1}^{ -1}  \tau^{2+\alpha+\beta}.									 		
				\end{aligned}
			\end{equation*}
			At this point, we employed \eqref{parabolicsmoothing1} and \eqref{eq71}. 

			After treating the boundary terms separately, the sum of the remaining terms with $j \neq 0$ and $j \neq n-1$ can be bounded using \eqref{parabolicsmoothing1}, \eqref{eq68a}, \eqref{eq68b}, \eqref{eq73}, and \eqref{eq72} as follows
			\begin{equation*}
				\begin{aligned}
					\left \| \sum_{j=1}^{n-2} \mathrm{e}^{-j \tau A}  \mathcal{R}_{n-j}^{[1]}  \right \| 
					& \leq  \tau  \sum_{j=1}^{n-2} \left\| \mathrm{e}^{-j \tau A}  b_2(-\tau A)  B \right\| \times   \\
					& \qquad \int_0^{c_2 \tau}\int_0^\xi  \left \| \mathrm{e}^{-\left( c_2 \tau-\xi \right) A} B  \right\|  (\xi - \sigma)  \left \| u^{ \prime \prime}\left(t_{n-j-1}+\sigma\right) \right \|  \mathrm{d} \sigma   \mathrm{d} \xi  \\					
					& \leq 	C \tau  \sum_{j=1}^{n-2} t_{j}^{-\gamma} t_{n-j-1}^{\alpha-1} 	\int_0^{c_2 \tau} \int_0^\xi \left( c_2 \tau-\xi \right)^{-\gamma} (\xi - \sigma)  \mathrm{d} \sigma   \mathrm{d} \xi  \leq 	C t_n^{-\gamma+\alpha} \tau^{3-\gamma} ,
				\end{aligned}
			\end{equation*}
			{\allowdisplaybreaks
				\begin{align*}
					& \left \| \sum_{j=1}^{n-2} \mathrm{e}^{-j \tau A}  \mathcal{R}_{n-j}^{[5]}   \right \|  \leq  \sum_{j=1}^{n-2} \left\| \mathrm{e}^{-j \tau A} A^{1-\zeta} \right \| \times 
					\\ & \qquad \int_0^{ \tau} \int_0^\xi  \left\| \mathrm{e}^{-\left(  \tau-\xi \right) A} A^{\zeta} \right \|    \frac{ (\xi - \sigma)^2}{2} \left \| A^{-1}BA^\beta \right \|  \left \| (A-B)^{-\beta} u^{\prime \prime \prime}\left(t_{n-j-1}+\sigma\right) \right \|  \mathrm{d} \sigma   \mathrm{d} \xi  \\					 
					&  \qquad \qquad \qquad \quad \, \, \, \, \, \leq	C  \sum_{j=1}^{n-2} t_{j}^{\zeta-1} t_{n-j-1}^{\alpha+\beta-2}	\int_0^{ \tau} \int_0^\xi  \frac{ (\xi - \sigma)^2}{2} \left(  \tau-\xi \right)^{-\zeta} \mathrm{d} \sigma   \mathrm{d} \xi \\
					&   \qquad \qquad \qquad \quad \, \, \, \, \, \leq 	C \tau^{\frac{7}{2}-3\zeta} \sum_{j=1}^{n-2} t_{j}^{\zeta-1} t_{n-j-1}^{-1 +\zeta}    \leq 	C t_n^{-1+2\zeta } \tau^{\frac{5}{2}-3\zeta}  .
				\end{align*}
			}Here, we assume that there exist some $0< \alpha \leq \frac{1}{2}$ and $ \beta \in \mathbb{R}$ for which $\alpha+\beta = \frac{1}{2}-\zeta$. The exact values of $\alpha,\beta$ will be determined later. We are now in a position to state the main result.
			\begin{thm}\label{theo4} 
				Let Assumptions \ref{ass1} and \ref{ass2}  be satisfied. Assume that there exist some \( \beta, \widetilde{\Gamma} \in \mathbb{R} \) such that \( A^{-1} B A^{\frac{1}{2} \widetilde{\Gamma}} \), \( A^{\frac{1}{2} \widetilde{\Gamma}} B A^{-1} \), and \( A^{-1} B A^\beta \) are bounded. Additionally, suppose that there exist some \( 0 < \alpha \leq \frac{1}{2} \) such that \( u_0 \in \mathcal{D}(A^{1+\alpha}) \) and \( \alpha + \beta = \frac{1}{2}^- \). Assume that the order conditions No. 1-4 in \eqref{Table1} are satisfied in the strong form, whereas No. 5 in \eqref{Table1} is satisfied in the very weak form $A=0$. Then the numerical solution of the initial value problem \eqref{4121} using third-order exponential Runge–Kutta methods with $s=3$ in \eqref{eq20} satisfies the error bound
				\[
				\left\|u_n - u\left(t_n\right)\right\| \leq C t_n^{-1^+} \tau^{\frac{5}{2}^-}
				\]
				for $0 < \gamma \leq \frac{1}{2}$. For $\frac{1}{2} < \gamma \leq 1$, we have
				\[
				\left\|u_n - u\left(t_n\right)\right\| \leq C t_n^{-\gamma + \alpha} \tau^{3 - \gamma}
				\]
				uniformly for $0 \leq n\tau \leq T$ in both cases. The constant $C$ depends on $T$, but it is independent of $n$ and $\tau$.
			\end{thm}
			\begin{proof} First we consider the case $0<\gamma \leq \frac{1}{2}$. 
				It is obvious that: $ \forall 1 \leq j \leq n-2:   t_{n-j-1}^{-\gamma} \geq T^{-\gamma}:=\widetilde{C}.$ So, we get
				$ 1 = \frac{\widetilde{C}}{\widetilde{C}} \leq \frac{1}{\widetilde{C}} t_{n-j-1}^{-\gamma}. $
				By applying the norm in \eqref{eq18} and using the triangle inequality, the estimatation of $e_n$ has the following form
				\begin{equation*}
					\begin{aligned}
						\left \| e_n \right \| 
						& \leq  C \tau \sum_{j=1}^{n-2} t_{n-j-1}^{-\gamma}   
						\left \| e_j \right \|  + C \tau^{1-\gamma}     \left \| e_{n-1} \right \|  + \frac{C}{\widetilde{C}} \tau  \sum_{j=1}^{n-1}  t_{n-j-1}^{-\gamma}  \left \| e_j \right \|    + C t_n^{-1^+ } \tau^{\frac{5}{2}^-}  \\					
						& \leq  \widehat{C} \tau \sum_{j=1}^{n-2} t_{n-j}^{-\gamma}   
						\left \| e_j \right \|  + \widehat{C} \tau^{1-\gamma}     \left \| e_{n-1} \right \|   + \widehat{C} t_n^{-1^+ } \tau^{\frac{5}{2}^-}.			
					\end{aligned}
				\end{equation*}
			For the case $\frac{1}{2} < \gamma \leq 1$, we follow the same procedure. The only change in the proof is that the last term becomes $C t_n^{-\gamma + \alpha} \tau^{3 - \gamma}$ instead of $C t_n^{-1^+} \tau^{\frac{5}{2}^-}$.
				The proof is completed by applying a discrete Gronwall lemma (see \cite{HO2010}).
			\end{proof}
			In the finite-dimensional case, where \( A \) is a second-order differential operator with homogeneous Dirichlet boundary conditions and \( B \) is a first-order differential operator, Theorem \ref{theo4} can be applied with \( X = L^2(\Omega) \), where \( \alpha = \frac{1}{4} - \zeta \) (see \cite{HOCHBRUCK2005323}) and \( \beta = \frac{1}{4} - \zeta \) as established by Lemma \ref{Fourierlemma}. Therefore, we have $ \alpha + \beta = \frac{1}{2}^- $. Note also that $ A^{\frac{1}{2} \widetilde{\Gamma} }  B A^{-1} $ and $ A^{-1} B A^{\frac{1}{2}\widetilde{\Gamma}}  $, where $ \widetilde{\Gamma} = \frac{1}{2} - \zeta $, are bounded by Lemma \ref{Fourierlemma}, which verifies the estimations in Sections IV and V.

To illustrate the order of the bound given in Theorem \ref{theo4}, numerical investigations will also be conducted in spaces other than \( X = L^2(\Omega) \), including \( X = L^1(\Omega) \) and \( X = L^{\infty}(\Omega) \). The numerical results will show perfect agreement with the error bounds derived in Theorem \ref{theo4}.
		\begin{remark}
			One might ask whether it is possible to achieve the full order of third-order method by imposing a stronger assumption on the initial data. For example, we assume that \( u(0,x) \in \mathcal{D}(A^2)\cap \mathcal{D}((A-B)^2) \). We observe that the second-order derivative of \( u \) can be estimated using \eqref{parabolicsmoothing1} as follows
			\begin{equation}\label{eq12} 
			\left \| u^{\prime \prime }(\sigma) \right \| \leq \left \|  \mathrm{e}^{-\sigma(A-B)} (A-B)^{2}u_0\right \| \leq \left \| \mathrm{e}^{-\sigma(A-B)} \right \| \left \| (A-B)^{2}u_0\right \| \leq C.	
			\end{equation}
			As a result of \eqref{eq12}, the order of certain terms can be improved. For instance, the term $\mathrm{e}^{-(n-1) \tau A} \mathcal{R}_{1}^{[5]}$ can be estimated by employing \eqref{parabolicsmoothing1} and \eqref{eq15}  
			\begin{equation*}
				\begin{aligned}
					\left \|  \mathrm{e}^{-(n-1) \tau A} \mathcal{R}_{1}^{[5]}  \right \| 
					& \leq  \left \| \mathrm{e}^{-(n-1) \tau A} A^{\gamma}  \right \|  \int_0^{ \tau} \int_0^\xi \left\| \mathrm{e}^{-\left(  \tau-\xi \right) A} \right\|    (\xi - \sigma) \left \| A^{-\gamma} B \right \| \left \| u^{ \prime \prime}\left(\sigma\right) \right \|  \mathrm{d} \sigma   \mathrm{d} \xi  \\ &  \leq 	C  t_{n-1}^{ -\gamma} \tau^{3}.	
				\end{aligned}
			\end{equation*} 
However, some terms, such as IV and V, remain unaffected by the increased regularity of the initial data. Therefore, although the stronger assumption on \( u_0 \) may address certain terms, it may not resolve all of them.
		\end{remark}
		\begin{remark}
We can derive a method to satisfy conditions No. 1-4 in \eqref{Table1} in a strong form, and condition No. 5 in \eqref{Table1} in a weak form, where \( A = 0 \) is evaluated only at \( b_i(-\tau A) \), $2\leq i \leq s$. It is straightforward to verify that the following scheme aligns with these requirements
		{			
		\begin{equation}\label{eq23}
			\begin{array}{c|ccc}
				0 & & & \\
				\frac{1}{2} & \frac{1}{2} \varphi_{1}(-\frac{1}{2} \cdot) & & \\
				1 & \varphi_{1} - 2\varphi_{2}(-\frac{1}{2} \cdot) - 2\varphi_{2} & 2\varphi_{2}(-\frac{1}{2}  \cdot) + 2\varphi_{2} & \\
				\hline & \varphi_{1} - 3\varphi_{2} + 4\varphi_{3} &  4\varphi_{2} - 8\varphi_{3} & - \varphi_{2} + 4\varphi_{3}
			\end{array} ,
		\end{equation}
		}where $\varphi_{j}(-c \cdot)=\varphi_j\left(-c\tau A\right)$ and $ \varphi_j=\varphi_j(-\tau A)$.  For \( X = L^2(\Omega) \), the following estimate can be improved by utilizing \eqref{parabolicsmoothing1}, \eqref{eq15} (with \( \gamma = \frac{1}{2} \)), \eqref{eq73}, and Lemma \ref{Fourierlemma}.
		\begin{equation*}
			\begin{aligned}
				& \left \| \sum_{j=1}^{n-2} \mathrm{e}^{-j \tau A}  \mathcal{R}_{n-j}^{[1]}  \right \| \leq  C \tau  \sum_{j=1}^{n-2} \left\| \mathrm{e}^{-j \tau A} A      \right\| \left \| A^{-1} B A^{\frac{1}{4}-\zeta} \right \| \times \\ & \qquad  \int_0^{c_2 \tau} \int_0^\xi  \left \| \mathrm{e}^{-\left( c_2 \tau-\xi \right) A} A^{\frac{1}{4}+\zeta}  \right\|  (\xi - \sigma) \left \| A^{-\frac{1}{2}} B \right \| \left \| u^{ \prime \prime}\left(t_{n-j-1}+\sigma\right) \right \|  \mathrm{d} \sigma   \mathrm{d} \xi  \\				
				& \qquad \qquad \qquad \quad \, \, \, \, \leq 	C \tau  \sum_{j=1}^{n-2} t_{j}^{-1} t_{n-j-1}^{\alpha-1} 	\int_0^{c_2 \tau} \int_0^\xi \left( c_2 \tau-\xi \right)^{-\frac{1}{4}-\zeta} (\xi - \sigma)  \mathrm{d} \sigma   \mathrm{d} \xi \\
				& \qquad \qquad \qquad \quad \, \, \, \, \leq 	C t_n^{-1+\alpha } \tau^{\frac{11}{4}-\zeta} .
			\end{aligned}
		\end{equation*}
Similarly, we also obtain a convergence rate of \( \frac{11}{4}^- \) for IV. As a result, we would expect to achieve an improved order of two and three-quarters. The numerical investigations to confirm this are provided, as illustrated in Figure \ref{fig8}. It is beyond the scope of this paper to re-examine all the estimates here in this norm, rather than just give an example where the bound can be sharper. 	 				
		\end{remark}		
			\begin{lem}\label{Fourierlemma}
Let \(\zeta > 0\) be a fixed small number, and let \(\Omega\) be a bounded domain in \(\mathbb{R}^N\), where \(N\) denotes the dimension. Suppose that \( A \) is a second-order differential operator with homogeneous Dirichlet boundary conditions, and \( B \) is a first-order differential operator. Then $ A^{-1}BA^\beta $ is bounded for $\beta = -\zeta$ in $L^1(\Omega)$, $L^{\infty}(\Omega)$, and for $\beta = \frac{1}{4}-\zeta$ in $L^2(\Omega)$.
			\end{lem}
			\begin{proof}
				First, we note that, based on the assumptions regarding the operators \( A \) and \( B \) in this lemma, they do not commute. Without loss of generality, we express the partial sum of the Fourier series of any function $f$ (possibly with an odd extension) as \eqref{eq22}
				\begin{equation}\label{eq22}
f(\mathbf{x})=\sum_{\mathbf{k} \in \mathcal{I}_{\mathbf{M}}} \hat{f}(\mathbf{k}) \prod_{i=1}^N \sin \left(k_i \pi x_i\right)
				\end{equation}
where $\mathcal{I}_{\mathbf{M}}=\left\{\mathbf{k}=\left(k_1, k_2, \ldots, k_N\right) \mid 1 \leq k_i \leq M_i, i=1,2, \ldots, N\right\} 
$, $\mathbf{k} \in \mathcal{I}_{\mathbf{M}}$ indicates that $\mathbf{k}=\left(k_1, k_2, \ldots, k_N\right)$ satisfies $1 \leq k_i \leq M_i$ for all $i$, and  $\mathbf{x}=\left(x_1, x_2, \ldots, x_N\right) \in \Omega \subset \mathbb{R}^N$. 						
				Applying the operator $A^{-1}BA^\beta$ to \eqref{eq22} yields
				\begin{equation}\label{eq56}
	\begin{aligned}
					A^{-1}BA^\beta f(\mathbf{x}) &= \sum_{\mathbf{k} \in \mathcal{I}_{\mathbf{M}}} \hat{f}(\mathbf{k}) (k_i\pi)^{2\beta-1} \times \\ & \prod_{\substack{j=1 \\ j \neq i}}^N \left(  \sin \left( k_j \pi x_j \right) \cos \left( k_i \pi x_i \right) + (1-(-1)^k)x_i\sin \left( k_j \pi x_j \right) - \sin \left( k_j \pi x_j \right)   \right).	
\end{aligned}					
				\end{equation}
In $L^1(\Omega)$ norm, according to the Riemann–Lebesgue lemma (see \cite{1684b3a8-48ea-3eb7-8a1e-11227e20bf15}), we obtain
				\begin{align*}
					\left \| A^{-1}BA^{-\zeta}f(\mathbf{x}) \right \|_1 \leq C \int_{\Omega}   \sum_{\mathbf{k} \in \mathcal{I}_{\mathbf{M}}}  \left| \hat{f}(\mathbf{k}) \right|  (k_i\pi)^{-2\zeta-1}  \mathrm{d}\mathbf{x} \leq  C \pi^{-2\zeta-1} \sum_{\mathbf{k} \in \mathcal{I}_{\mathbf{M}}} k_i^{-2\zeta-1}  .
				\end{align*}
				For $L^2(\Omega)$ norm, using the fact that $ \hat{f}(\mathbf{k})$ is in  $l^2(\Omega)$ (see \cite{garfken67:math}), we get 
				\begin{align*}
					\left \| A^{-1}BA^{\frac{1}{4}-\zeta}f(\mathbf{x}) \right \|_2^2  \leq C \int_{\Omega}  \sum_{\mathbf{k} \in \mathcal{I}_{\mathbf{M}}}  \hat{f}(\mathbf{k})^2   \sum_{\mathbf{k} \in \mathcal{I}_{\mathbf{M}}}  (k_i\pi)^{-4\zeta-1}   \mathrm{d}\mathbf{x} \leq  C \sum_{\mathbf{k} \in \mathcal{I}_{\mathbf{M}}}  \hat{f}(\mathbf{k})^2   \sum_{\mathbf{k} \in \mathcal{I}_{\mathbf{M}}}  (k_i\pi)^{-4\zeta-1}  .
				\end{align*}
				Lastly, for estimation in  $L^{\infty}(\Omega)$ norm, by the definition of the Fourier coefficients, we have $\left\|f_k\right\|_{\infty} \leq C.$ An upper bound for \( A^{-1} B A^{-\zeta} f(\mathbf{x}) \) is established as follows
				\begin{align*}
					\left \| A^{-1} B A^{-\zeta} f(\mathbf{x}) \right \| \leq C   \sum_{\mathbf{k} \in \mathcal{I}_{\mathbf{M}}}  \left| \hat{f}(\mathbf{k}) \right|  (k_i\pi)^{-2\zeta-1}   \leq  C   \sum_{\mathbf{k} \in \mathcal{I}_{\mathbf{M}}}    (k_i\pi)^{-2\zeta-1} 
				\end{align*}
As \( M_1,\ldots,M_N \to \infty \), the desired results for various norms are obtained. 
			\end{proof}
			\section{Numercal investigations}\label{numericalchap2}
			This section presents numerical experiments to verify the error bounds derived. First, we  discuss the implementation of various exponential Runge--Kutta schemes.
			\subsection{Implementation}
			Implementing exponential Runge--Kutta methods requires approximating the action of a matrix function on a vector. This involves the linear combinations of several \( \varphi_i(-\tau A) \) functions with vectors, represented as \( \sum_{i=1}^k \varphi_i(-\tau A) v_i \). We transform \( \sum_{i=1}^k \varphi_i(-\tau A) v_i \) into a single matrix exponential applied to a vector using an augmented matrix (see \cite{doi:10.1137/100788860}), denoted as \( \exp(-\tau \widehat{A}) V_0 \). 
			

			To implement the second-order exponential Runge--Kutta methods (see \cite{HO2005}), we set the parameter \( c = \frac{1}{2} \), while for the third-order methods we use the scheme \(\mathtt{ETD3RK}\) (see \cite{HO2005}). The reference solution is computed using the RK4 method with a sufficiently small time step.

			\subsection{Numerical results}
Throughout the paper, the proof holds for any finite dimension. Therefore, it suffices to perform numerical experiments with the following 1D linear advection-diffusion problem
			\begin{equation}\label{problem2}
				\begin{aligned}
					& \partial_t  u(t,x) - 0.2 \partial_{xx} u (t,x)  =\partial_x u(t,x) ,\qquad (t,x)\in[0,1]\times[0,1], 
				\end{aligned}	
			\end{equation} 
			subject to homogeneous Dirichlet boundary conditions. We employ the standard second-order finite difference scheme to discretize the diffusion part $A= -0.2 \partial_{xx}$ and the advection part $B=\partial_{x}$, with $n=399$ inner grid points, resulting again in matrices $A$ and $B$, respectively. Note that the operator \( A \) includes  the homogeneous Dirichlet boundary conditions; therefore, the operators \( A \) and \( B \) do not commute in both the continuous and discrete cases.
			 We have \( A^{-\frac{1}{2}}B \) and \( BA^{-\frac{1}{2}} \) are bounded operators that satisfy Assumption \ref{ass2} with \( \gamma = \frac{1}{2} \). Therefore, based on Theorems \ref{theo2} and \ref{theo3}, we get the original order for the first and second order methods. In addition, according to Theorem \ref{theo4}, we would expect a reduction to two and a half for the third order method. 
			
\begin{figure}[h!]
				\centering
				\subfigure[\label{fig7a}]{\includegraphics[width=0.48\textwidth]{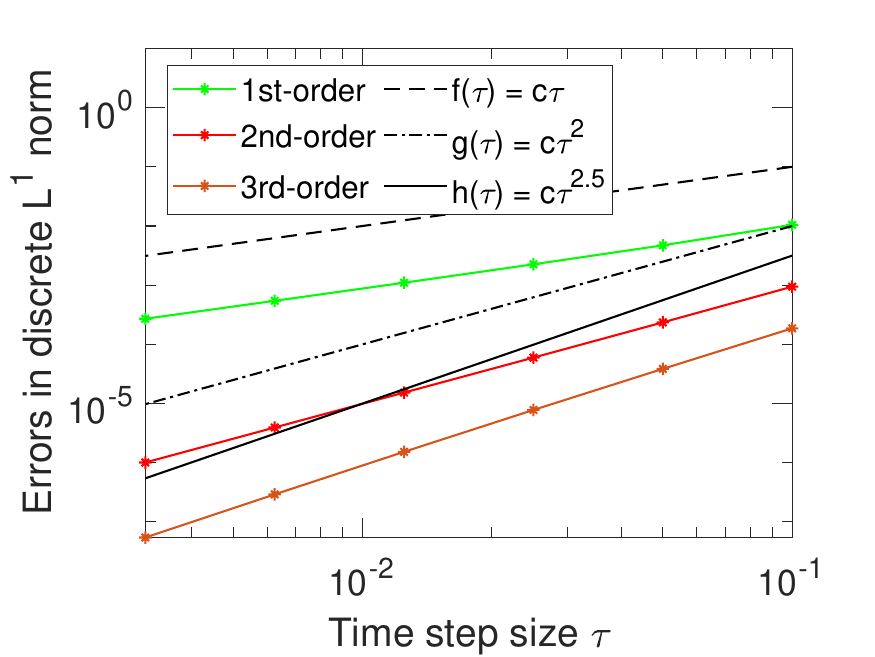}} 
				\subfigure[\label{fig7b}]{\includegraphics[width=0.48\textwidth]{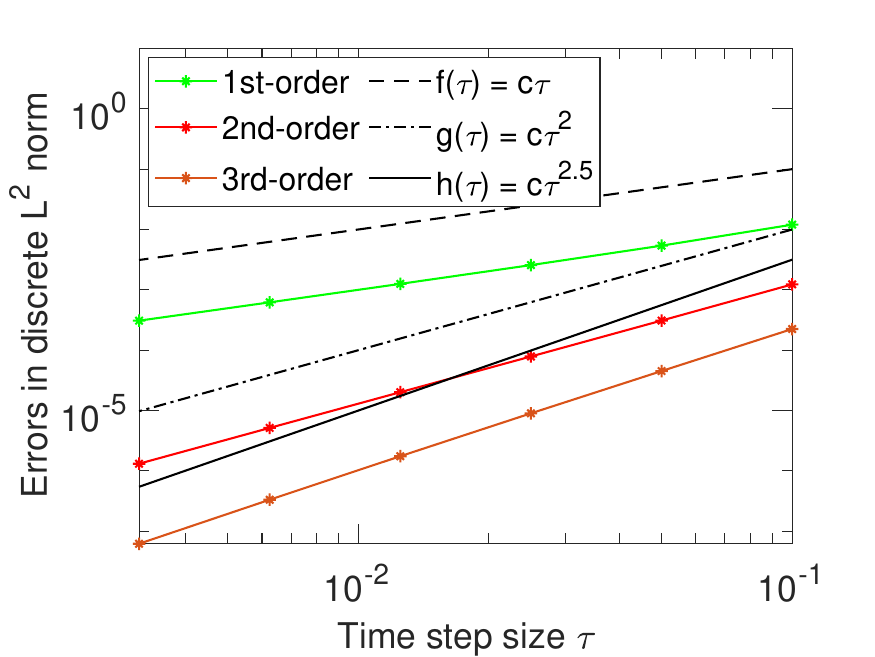}}	\\		
				\centering
				\subfigure[\label{fig7c}]{\includegraphics[width=0.48\textwidth]{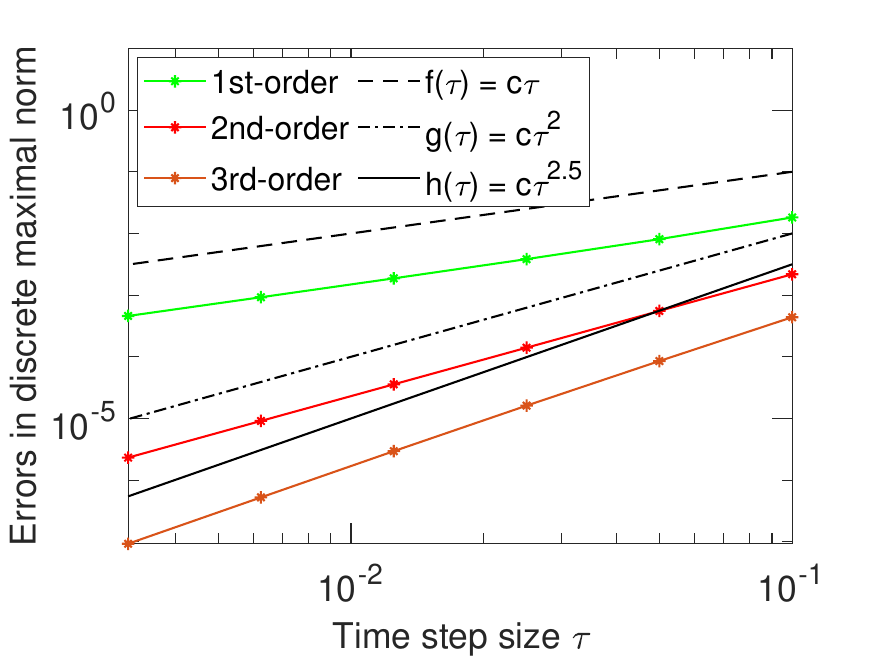}} 
				\caption{ The global error is shown as a function of the time step \( \tau \) for the different methods. '1st-order' refers to the first-order method, with similar interpretations for higher-order methods.}
				\label{fig7} 				
			\end{figure}

\begin{figure}[h!]
	\centering
	\includegraphics[width=0.7\textwidth]{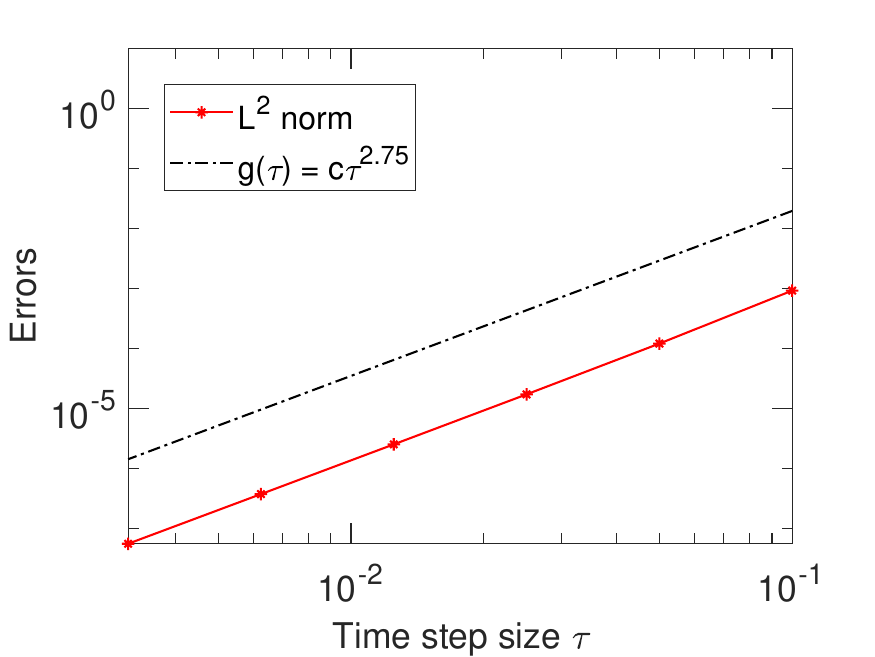}
	\caption{ The global error is shown as a function of the time step \( \tau \) for the method given in \eqref{eq23}. }
	\label{fig8}
\end{figure}

			We choose the initial data $4x(1-x) = u(0,x) \in \mathcal{D}(A)$. With this initial data, it also satisfies part of Theorems \ref{theo2}, \ref{theo3}, and \ref{theo4}.
			 The plots illustrating the numerical results for different norms are shown in Figures \ref{fig7a}, \ref{fig7b} and \ref{fig7c}. These plots are in complete agreement with Theorems \ref{theo2}, \ref{theo3}, and \ref{theo4}. The numerical results in higher dimensions are similar; therefore, it does not make sense to perform additional numerical experiments.
			
			\section{Conclusion}\label{concluchap2}
			In this paper, we investigate the application of exponential Runge--Kutta methods to the initial value problem \eqref{4121} involving unbounded operators and non-commuting operators \( A \) and \( B \). We treat the dominant part $A$ of \eqref{4121} exactly, wheareas handling the unbounded operator $B$ explicitly. We present a careful error analysis for the application of these methods up to third-order to the initial value problem \eqref{4121} within the framework of analytic semigroups in an abstract Banach space. For the first- and second-order exponential Runge–Kutta methods, no order reduction takes place. In the convergence analysis of the first-order method, we present an additional result without the assumption that the initial data belongs to the domain of \( A \).
			For third-order methods, the convergence analysis reveals where the order reduction happens. In an abstract Banach space \( X \), where the third-order method satisfies the last order condition in a very weak form, we  achieve an order of two and a half. However, if \( X = L^2(\Omega) \) and the newly presented scheme is used, we obtain a slightly higher order of two and three-quarters. To validate sharpness of convergence bounds, numerical investigations are conducted. Overall, numerical results are in excellent agreement with theoretical findings. 

\section*{Acknowledgements}

The author is deeply grateful to Professor Alexander Ostermann and Professor Lukas Einkemmer for their unwavering guidance, and numerous constructive discussions.

\bibliography{references}

\end{document}